\documentclass[11pt]{article}
\usepackage[latin1]{inputenc}
\usepackage{amsmath,amssymb}
 \usepackage{xcolor}
\usepackage{multicol}
\usepackage{floatrow}
\usepackage{latexsym,epsfig}
\usepackage{amssymb,amsmath,amsfonts,amscd}
\usepackage{exscale}
\usepackage{ifthen}
\usepackage{longtable}
\usepackage{subfigure}
 \usepackage{cite}
\usepackage{lscape}
\usepackage{blindtext}
\usepackage{soul}
\usepackage{mathtools}

\usepackage{ulem}

 \newtheorem{theo}{\textbf{Theorem}\ }
[section]
\newtheorem{lem}[theo]{\textbf{Lemma}\ }
\newtheorem{coro}[theo]{Corollary\ }

\newtheorem{prop}[theo]{\textbf{Proposition}\ }

\newtheorem{notation}[theo]{Notation\ }

\newtheorem{property}[theo]{\textbf{Property}\ }

\def \N{\mathbb{N}}

\def \E{\mathbb{E}}
\def \X{\mathbb{X}}
\def \R{\mathbb{R}}

\def \P{\mathbb{P}}
\def \S{\mathcal{S}}

	\definecolor{amber}{rgb}{1.0, 0.75, 0.0}
	\definecolor{airforceblue}{rgb}{0.36, 0.54, 0.66}
		\definecolor{alizarin}{rgb}{0.82, 0.1, 0.26}
	\definecolor{amaranth}{rgb}{0.9, 0.17, 0.31}	
	\definecolor{amber(sae/ece)}{rgb}{1.0, 0.49, 0.0}	
	\definecolor{americanrose}{rgb}{1.0, 0.01, 0.24}	
	\definecolor{amethyst}{rgb}{0.6, 0.4, 0.8}	
	\definecolor{antiquebrass}{rgb}{0.8, 0.58, 0.46}	
	\definecolor{ao(english)}{rgb}{0.0, 0.5, 0.0}	
	\definecolor{apricot}{rgb}{0.98, 0.81, 0.69}	
	\definecolor{aquamarine}{rgb}{0.5, 1.0, 0.83}	
	\definecolor{armygreen}{rgb}{0.29, 0.33, 0.13}
	\definecolor{arsenic}{rgb}{0.23, 0.27, 0.29}	
	\definecolor{asparagus}{rgb}{0.53, 0.66, 0.42}	
	\definecolor{auburn}{rgb}{0.43, 0.21, 0.1}	
	\definecolor{awesome}{rgb}{1.0, 0.13, 0.32}	
	\definecolor{ballblue}{rgb}{0.13, 0.67, 0.8}	
	\definecolor{bananayellow}{rgb}{1.0, 0.88, 0.21}	
	\definecolor{bittersweet}{rgb}{1.0, 0.44, 0.37}	
	\definecolor{bleudefrance}{rgb}{0.19, 0.55, 0.91}
	\definecolor{blush}{rgb}{0.87, 0.36, 0.51}	
	\definecolor{bostonuniversityred}{rgb}{0.8, 0.0, 0.0}	
	\definecolor{brickred}{rgb}{0.8, 0.25, 0.33}	
\definecolor{brightmaroon}{rgb}{0.76, 0.13, 0.28}	
	\definecolor{brilliantrose}{rgb}{1.0, 0.33, 0.64}
	\definecolor{britishracinggreen}{rgb}{0.0, 0.26, 0.15}	
	\definecolor{bulgarianrose}{rgb}{0.28, 0.02, 0.03}	
	\definecolor{burgundy}{rgb}{0.5, 0.0, 0.13}	
	\definecolor{burntorange}{rgb}{0.8, 0.33, 0.0}	
	\definecolor{byzantium}{rgb}{0.44, 0.16, 0.39}			
	\definecolor{candyapplered}{rgb}{1.0, 0.03, 0.0}
	\definecolor{capri}{rgb}{0.0, 0.75, 1.0}	
	\definecolor{carminered}{rgb}{1.0, 0.0, 0.22}	
	\definecolor{charcoal}{rgb}{0.21, 0.27, 0.31}	
	\definecolor{coquelicot}{rgb}{1.0, 0.22, 0.0}	
	\definecolor{daffodil}{rgb}{1.0, 1.0, 0.19}

\numberwithin{equation}{section}

\begin{document}
\title{Central limit theorem  
	for a critical multi-type  branching process  
in random environment }

  \maketitle 

\begin{center}
	{\Huge \bf }  
 E. Le Page $^($\footnote{
	Universit\'{e} de Bretagne-Sud, LMBA UMR CNRS 6205, Vannes, France. \\ emile.le-page@univ-ubs.fr}$^)$,
M. Peign\'e  
 $\&$  C. Pham 
$^($\footnote{
Institut Denis Poisson UMR 7013,  Universit\'e de Tours, Universit\'e d'Orl\'eans, CNRS  France. \\ 	 peigne@univ-tours.fr, thidacam.pham@univ-tours.fr}$^)$

\vspace{0.5cm}
     \end{center}

         \centerline{\bf  \small Abstract }

  Let  $(Z_n )_ {n \geq 0} $  with $Z_n= (Z_n(i, j))_{1\leq i, j \leq p}$ be a $p$ multi-type critical branching process in random environment, and let $M_n$ be the   expectation of  $Z_n$ given a fixed environment. We prove   theorems on convergence  in distribution   of sequences of branching processes $\left\{ {Z_{n}\over \vert M_{n}\vert} \slash \vert Z_{n}\vert >0\right\}$ and 
  $\left\{ {\ln Z_{n}\over \sqrt{n}} \slash \vert Z_{n}\vert >0\right\}$. These theorems extend  similar results  for single-type critical branching process in random environment.

  \vspace{0.5cm}

  \noindent Keywords: multi-type branching process, random environments, central limit theorem.


\section{Introduction}

Single-type and multi-type branching processes in random environments  (BPREs) are a central topic of  research; they were introduced in the $1960$s  in order to  describe the development of populations whose evolution may be randomly affected by environmental factors.

 In the single-type case, the behaviour of these processes is mainly determined  by the $1$-dimensional random walk generated by the logarithms of the expected population sizes $m_k$, for $  k \geq 0,$ of  the respective generations; they are classified  in three classes - supercritical, critical and subcritical - of single-type BPREs,  according to the fact that the associated random walk tends to $+\infty$, oscillates or tends to $-\infty$.  Their study is closely related to the   theory of  fluctuations of random walks on $\mathbb R$ with i.i.d. increments;  when $\mathbb E\left[ \vert \log m_k \vert \right]<+\infty$, the BPRE is supercritical (resp.,  critical and subcritical)  when      $\mathbb E\big[ \vert  \log m_k \vert  \big]>0$  (resp.,  $\mathbb E \big[\vert \log m_k \vert \big]=0$ or $\mathbb E \big[\vert \log m_k \vert \big]<0$).

In this context, a huge body of papers is devoted to study of   the asymptotic behaviour  of the probability of non-extinction up to time $n$ and the  distribution of the population size  conditioned to survival up to time $n$. In the critical case,   the branching process is degenerate with probability one  and the   probability of non-extinction up to time $n$  is equivalent to $c /\sqrt{n}$  as $n\to +\infty$, for some explicit positive constant $c $ \cite{Ko}, \cite{GK}. The convergence in distribution of the process  conditioned to non-extinction comprises first  the Yaglom classical theorem; the convergence of finite-dimensional distributions of the processes  was established by Lamperti and Ney \cite{LN} who showed that the limiting process   is a diffusion process and described its transition function. 
V.I.  Afanasev   described  the limiting process  in terms of Brownian excursions  \cite{A}. These statements   more or less claim  that  the conditional logarithmic behaviour of the BPRE given its non-extinction at the terminal time $n$ is the same as the one of the associated random walk of conditional mean values conditioned to staying  positive. These results are extended in \cite{AGKV2005}  under more general assumptions, known as Spitzer's condition in fluctuation theory of random walks, and some additional moment conditions.

It is of interest to prove analogues of the above statements for   multi-type BPREs $(Z_n)_{n \geq 0}$. As in the single-type case, the set of multi-type BPREs may be divided  into three  classes: they are supercritical (resp.  critical or subcritical) when the upper Lyapunov exponent of the product of  random mean matrices $M_k$  is positive (resp. null or negative) \cite{Kaplan}. Let us emphasize that the role of the random walk associated to the BPRE  in the single-type case   is played  in the multi-type case by the logarithm of the norm of some $\mathbb R^p$-valued Markov chain whose increments  are governed by i.i.d. random $p\times p$-matrices $M_k$, for $ k \geq 0$, whose coefficients are non-negative and correspond to  the expected population sizes of the respective generations, according to the types of the particles and their direct parents.  Product of random matrices is the object of several investigations and many limit theorems do exist in this context  (see   \cite{BL}  and references therein). The theory of their  fluctuations   is recently studied during the last decade  using the promising approach initiated by V. Denisov   V. Wachtel \cite{DW}.

The question of the asymptotic behaviour of the  probability  of non-extinction up to time $n$ is solved recently, under quite general moment assumptions   and irreducibility condition  on the projective action  of the product of matrices $M_k$; as in the single-type case, it is proved that the   probability of non-extinction up to time $n$  is equivalent to $c /\sqrt{n}$  as $n\to +\infty$, for some explicit positive constant $c$  \cite {LPP}, \cite{DV}.   The asymptotic distribution of   the size of the population conditioned to non-extinction remains open; 
in this paper, we prove a central limit theorem for the logarithm of the size of the population at time $n$, conditioned to non-extinction.

\section{Preliminaries, hypotheses and statements}

We fix an integer $p\geq 2$ and  denote by $\mathbb R^p$ (resp. $\mathbb N^p$) the set of $p$-dimensional column vectors with real (resp.  non-negative integer) coordinates; for any column vector $x$ of $\R ^p$ defined by $x = (x_i)_{1\leq i\leq p}$, we denote by $\tilde x$ the row vector $\tilde x:= (x_1, \ldots, x_p)$. Let  $  {\bf 1}$ (resp. $  {\bf 0}$)  be the column vector of $\mathbb R^p$ whose all coordinates equal  $1$ (resp. $0$).
We denote by $ \left\{{e}_i, 1 \le i \le p\right\}$  the canonical basis, by $\langle \cdot, \cdot\rangle $ the usual scalar product
 and by $|.|$ the corresponding $ \mathbb L^1$ norm. We also consider the general linear  semi-group $\S^+$   of $ p \times p$ matrices with  non-negative coefficients. We endow  $\S^+$ with the $\mathbb L^1$-norm denoted also by $|.|$.

The multi-type Galton-Watson  process we study here  is the  Markov chain   $(Z_n)_{n \geq 0}$  whose states are $p\times p$ matrices with integer  entries. We always assume that $Z_0$ is non-random. For any $1\leq i, j\leq p$,   the $(i, j)^{th}$  component $Z_n (i, j)$ of $Z_n$ may be interpreted   as the number of particles of type $j $ in the $n^{th}$ generation  providing that the ancestor at time $0$ is of  type $i$. In particular, $  \sum_{i=1}^pZ_{n}(i, j)=\vert Z_n(., j)\vert $ equals the number  of particles  of type $j$  in the $n^{th}$ generation  when there is 1  ancestor  of each type  at generation $0$; this quantity equals $ Z_n(\tilde {\bf 1}, j) $, with the notations introduced below. Similarly, 
 $  \sum_{j=1}^pZ_{n}(i, j)= \vert Z_n(i, \cdot)\vert$ equals the size of the population at time $n$ when there is one ancestor  of type $i$ at time 0.

Let us introduce Galton-Watson process in varying environment and assume that the offspring distributions of the process   $(Z_n)_{n \geq 0}$ are given by a sequence of  $\mathbb N^p$-valued  random variables $(\xi_n)_{n \geq 0}$.   More precisely, the distribution of the
number of typed $j$ children born to a single-typed $i$ parent at time $n$ is the same as the one of  $\xi_n(i, j)$.  
Let $ \xi_{n, k}, n\geq 0,  k \geq 1,$ be independent random variables defined on $(\Omega, \mathcal F, \mathbb P)$, where  the $\xi_{n, k}, k \geq 1$ have the same distribution as $\xi_n$, for any $ n \geq 0$.

The process $(Z_n)_{n \geq 0}$ is thus defined by recurrence as follows: for any $1\leq i, j\leq p$ and $n \geq 0$,
\[
Z_{n+1}(i, j) := \left\{
\begin{array}{lll} 
\displaystyle  \sum_{\ell=1}^p\quad  \sum_{k= 1+ Z_{n}(i, 1)+\ldots +Z_{n}(i, \ell-1) }^{ Z_{n}(i, 1)+\ldots  +Z_{n}(i, \ell)}\quad \xi_{n, k}(\ell, j) & \text{when} &\displaystyle \sum_{\ell=1}^pZ_{n}(i, \ell)>0;
\\ \ &\ & \\
0& \text{otherwise.}& 
\end{array}
\right.
\]

We denote by $\mathcal G$ the set  of   multivariate  probability generating functions  $g =
(g^{(i)} )_{1\leq i \leq p}$ defined by: 
\[
g^{(i)} (s) = \sum_{ \alpha  \in {\mathbb{N}^p}} {{p^{(i)}  }(\alpha)}s^{ \alpha}, 
\]
for any $s=(s_i)_{1\leq i\leq p}\in [0, 1]^p$, where
\begin{enumerate}
	\item  $\alpha =(\alpha_i)_i  \in \mathbb{N}^p$ and $s^{  \alpha} = s_1^{{\alpha _1}} \ldots s_p^{{\alpha _p}}$;
	\item $p^{(i)}  (\alpha)=p^{(i)}  ({{\alpha _1}, \ldots,{\alpha _p}})$ is the  probability that a parent of type $i$  has $\alpha_1$ children of type $1,\ldots,\alpha_p$ children of type $p$.
\end{enumerate}
 
For each $1\leq i \leq p,$ the  distribution  of the $i^{\rm th}$ row vector  $(\xi_n(i, j))_{1\leq j\leq p}$ of the $\xi_n$ is characterized by  its generating function denoted by  $g_n^{(i)}$, and  set $g_n = (g_n^{(i)})_{1\leq i\leq p}$ and ${\bf g}=(g_n)_{n \geq 0}$. For any  $s=(s_i)_{1\leq i\leq p}\in [0, 1]^p$,
\[
g_n^{(i)}(s) = \mathbb E\left[ s_1^{\xi_n(i, 1)} s_2^{\xi_n(i, 2)} \ldots s_p^{\xi_n(i, p)}\right].
\]
For a given sequence ${\bf g}$ in $\mathcal G$,   we  denote by $Z^{\bf g}$ the Galton-Watson process corresponding to  ${\bf g}$ and omit the exponent ${\bf g}$ when there is no confusion; furthermore, we set  for any $n \geq 1$,    
\[g_{0,n}  = g_0\circ g_1  \ldots  \circ g_{n-1}=(g^{(i)}_{0,n})_{1\leq i\leq p},
\]
where $g^{(i)}_{0,n}=g_0^{(i)}\circ g_1  \ldots  \circ g_{n-1}$.
For any $1\leq i \leq p$ and  $ s\in [0, 1]^p$,
\begin{eqnarray*} \label{eqn1.1}
	\mathbb{E} \left[ s^{Z_n^{\bf g}(i, \cdot)}\Big\slash Z^{\bf g}_0(i, \cdot) \ldots, Z^{\bf g}_{n-1}(i, \cdot)\right] = g_{n-1}(s) ^{Z^{\bf g}_{n-1}(i, \cdot)} , 
\end{eqnarray*} 
which yields that
\[
	\mathbb{E} \left[s^{Z^{\bf g}_n(i, \cdot)} \right]=\mathbb E \left[s_1^{Z^{\bf g}_n(i, 1)} s_2^{Z^{\bf g}_n(i, 2)} \ldots s_p^{Z^{\bf g}_n(i, p)}\right] =  g_{0, n}^{(i)}(s) .
\]
More generally, for any $z = (z_1, \ldots, z_p) \in \mathbb N^p\setminus \{\tilde {\bf 0}\}$, we denote by $Z_n^{\bf g}(\tilde z, j)$ the number of particles of type $j $ in the $n^{th}$ generation   providing that there are $z_i$ ancestors of type $i$   at time $0$, for any $1\leq i \leq p$. 
Therefore, 
\begin{eqnarray*}
	\mathbb{E} \left[s^{Z^{\bf g}_n(\tilde z, \cdot)} \Big\slash Z^{\bf g}_0(\tilde z, \cdot), \ldots, Z^{\bf g}_{n-1}(\tilde z, \cdot)\right] = g_{0, n}(s) ^{ \tilde z}.
\end{eqnarray*}
Before going further, we introduce some necessary notations. For $n \ge 0$, 
\begin{enumerate}
	\item we denote by  $  M_{g_n}$ the mean matrix $\E \xi_n$:	
\[
M_{g_n} = (M_{g_n}(i, j))_{1\leq i, j \leq p} \quad \text{with} \quad M_{g_n}(i,j) = \E[\xi_n(i,j)];
\]
	in other words, 
	\[  M_{g_n}  = \left({\begin{array}{*{20}{c}}
		{\displaystyle \frac{{\partial g^{(1)}_{n} (\tilde {\bf 1})}}{{\partial {s_1}}}}& \ldots &{\displaystyle \frac{{\partial g^{(1)}_{n} (\tilde {\bf 1})}}{{\partial {s_p}}}}\\
		\vdots &{}&{}\\
		{\displaystyle \frac{{\partial g^{(p)}_{n} (\tilde {\bf 1})}}{{\partial {s_1}}}}& \ldots &{\displaystyle \frac{{\partial g^{(p)}_{n} (\tilde {\bf 1})}}{{\partial {s_p}}}}
		\end{array}} \right);\]

			\item for $1\leq i\leq p$,  let $ B_{g_n}^{(i)}$ be the Hessian matrices  
			\begin{align*}
	B_{g_n}^{(i)} &= \Bigl(B_{g_n}^{(i)}(k, \ell)\Bigr)_{1\leq k, \ell \leq p}:= \left({\partial^2g^{(i)}_{n} \over \partial s_k\partial s_\ell}(\tilde{\mathbf  1}) \right)_{1\leq k, \ell\leq p}  \\
	& =\Bigl(\mathbb E \big[\xi_n(i, k)(\xi_n(i, \ell) -\delta_{k, \ell})\big]\Bigr)_{1\leq k, \ell \leq p};
	\end{align*}	in particular, $\sigma^2_{g_n}(i, j):= {\rm Var} (\xi_n(i, j))=B_{g_n}^{(i)}(j, j)+M_{g_n}(i, j)-M_{g_n}(i, j)^2$.
\item and we set 	$\displaystyle \mu_{g_n} := \sum_{i=1}^p \bigl\vert B^{(i)}_{g_n} \bigr\vert$ and $\displaystyle \eta_{g_n}:= {\mu_{g_n}\over \vert M_{g_n}\vert ^2}$.
\end{enumerate}
 	The product of matrices  $M^{{\bf g}}_{0, n}:= M_{g_0}\ldots M_{g_{n-1}}$  controls the mean value of $Z_n$, according to the value of $Z_0$; indeed $\mathbb E[Z_n]= Z_0M_{g_0}\ldots M_{g_{n-1}} $ for any $n \geq 0$.
	The matrices $M_{g_n}$, for $ n\geq 0$, have non-negative entries which plays an important role on  the asymptotic  behaviour of the products $M^{{\bf g}}_{0, n}$, for $ n \geq 0$.
		
	We consider the cone   of $p$-dimensional row vectors
	\[
	\mathbb R^p_+ := \left\{ {\tilde x=(x_1,\ldots, x_p) \in \mathbb{R}^p \ \Big\slash \ \,\forall i = 1,\ldots,p,\,\,{x_i} \ge 0} \right\},
	\]
	and   the corresponding simplex  $\mathbb{X}$ defined by:
	\[
	{\mathbb X} := \{\tilde  x \in \mathbb R^p_+ \ \Big\slash\  \vert \tilde x\vert =1\}.
	\] 
 We introduce the actions of the semi-group $\S^+$ on $\mathbb R^p_+ $ and $\mathbb X$ as follows:
	\begin{itemize}
		\item the right and the left   linear actions  of    $\S^+$   on $\mathbb R^p_+ $ defined by: 
		\[
		(\tilde x,   M) \mapsto  \tilde x M\quad {\rm and} \quad (\tilde x,   M) \mapsto Mx 
		\]
		for any $\tilde  x \in\mathbb R^p_+$ and $M\in \S^+$, 
		
		\item the right and the left   projective actions  of    $\S^+$   on $\mathbb X$ defined by: 
		\[
		(\tilde x,   M) \mapsto  \tilde x \cdot M := \displaystyle \frac{\tilde x M}{\vert\tilde x M\vert} \quad {\rm and} \quad (\tilde x,   M) \mapsto M\cdot x := \displaystyle \frac{  Mx}{\vert Mx\vert} 
		\]
		for any $\tilde  x \in  \mathbb{X}$ and $M \in \S^+$.
	\end{itemize}
	
	For any  $M \in \S^+$,    let   
	$\displaystyle v(M) := \min_{1\leq j\leq p}\Bigl(\sum_{i=1}^d M(i, j)\Bigr).
	$
	Then  for any $\tilde x \in \mathbb R^p_+$,
	\begin{equation*} 
	0< v(M)\  \vert x\vert \leq \vert Mx\vert \leq \vert M\vert \ \vert x\vert.
	\end{equation*}
	We set $\mathfrak n(M):= \max \left({1\over v(M)}, \vert M\vert\right) $.
	
	We also introduce some proper subset of $\S^+$ which is of interest in the sequel: for any  constant $B\geq 1$, let  $\S^+(B)$ denote  the set of $p\times p$ matrices   $M$ with positive coefficients   such  that   for  any  $1\leq  i, j, k, l \leq p$,   
	\begin{equation*} 
	\frac{1}{B} \leq \frac{ M(i,j)} {M(k,l)}  \leq B.
	\end{equation*}
	
 Following  \cite{DHKP}, we introduce  some proper subset  of  generating functions  of   offspring distributions;  let $\xi$ be a $\mathbb N^p$-valued random variable defined by $\xi =(\xi (i, j))_{1\leq i, j \leq p}$, with generating function $g $ (as described above).
	
  \begin{notation}  Let $\varepsilon \in ]0, 1[$ and $K>0$. We denote  by $\mathcal  G_{ \varepsilon, K}$ the set of  generating functions   of multivariate offspring distributions    satisfying the following non-degeneracy assumptions:     for any $n \ge 0$ and   $1\leq i, j \leq p,$
	\begin{enumerate}
		\item[(1)]$ \quad  \mathbb P \big(\xi (i, j)\geq 2 \big) \geq \varepsilon $,
		\item[(2)] $ \quad  \mathbb P\bigl( \xi (i, \cdot)  = {\bf 0} \bigr)\geq \varepsilon,$
		\item[(3)] 
		$  \quad  \E\bigl[ \vert \xi (i ,\cdot)\vert ^2\bigr]\leq K<+\infty.$
	\end{enumerate}
	\end{notation} 
	D. Dolgopyat and co-authors  in  \cite{DHKP} proposed  a deep and useful description of the behaviour of the process $(Z_n^{\bf g})_{n \geq 0}$ when  all the generating functions $g_n$ of  the varying environment $\bf g$ belong to $\mathcal  G_{ \varepsilon, K}$. We present and extend  their results in section \ref{Ontheprobabilityofextinctioninvaryingenvironment}, they play a key role in controlling the fluctuations  of the Galton-Watson process in random environment, conditioned to non-extinction (see Corollary \ref{corotildeP}).

{\bf In random environment}, we consider   a sequence ${\bf f}=(f_n)_{n \geq 0}$   of  i.i.d. random variables defined on     $ ({\Omega,\mathcal{F},\mathbb{P}})$ and we  set   $f_{0,n}  = f_0\circ f_1  \ldots  \circ f_{n-1}$ for any $n \geq 1$;   as above, for any  $\tilde z = (z_1, \ldots, z_p) \in \mathbb N^p\setminus \{\tilde {\bf 0}\}$ and  $ s\in [0, 1]^p$,
\[
	\mathbb{E} \left[s^{Z_n(\tilde z, \cdot)}\Big\slash Z_0(\tilde z, \cdot) \ldots, Z_{n-1}(\tilde z, \cdot); f_0, \ldots, f_{n-1}\right] =  f_{0, n}(s) ^{Z_{n-1}(\tilde z, \cdot)}. 
\]

 For $n \geq 0$, the random matrices  $ M_{n} $ and  $ B_n^{(i)}$  are i.i.d. and of non-negative entries.  The common law of the $ M_{n} $ is  denoted by $\mu$. In order to simplify the notations, we set $M_n:= M_{f_n}, B_n^{(i)}=B_{f_n}^{(i)}$ and $\eta_n:= \eta_{f_n}$, with the notice that $\eta_n $ are non-negative  real valued i.i.d. random variables. Moreover, let $M_{0, n}$ and $M_{n, 0}$ denote the right and the left product of random matrices $ M_k$ for $k \geq 0$, respectively
$ M_{0, n} =  M_{0, n}^{\bf f} =M_0 M_1 \ldots  M_{n-1}\, \mbox{and} \,\,\, M_{n, 0} =  M_{n-1} \ldots   M_1  M_0$, with the convention that $M_{0, 0}=\text{I}$.
 Therefore,
\[
\mathbb{E} \left[Z_n\Big\slash f_0, f_1, \ldots, f_{n-1}\right]= Z_0 M_0 \ldots M_{n-1}=Z_0M_{0, n} \qquad  \mathbb P\text{-a.s.}.
\]

For any $1\leq i\leq p, $ the  probability of non-extinction  at generation $n$ given  the environment $f_0, f_1, \ldots  f_{n-1}$ is 
\begin{eqnarray*} 
	q_{n, i}^{\bf f}&:=& \mathbb{P} \Big( Z_n(i, \cdot) \ne \tilde {\bf 0}\Big\slash f_0 , f_1,\ldots, f_{n-1} \Big) \notag \\
	&=& 1-f_0^{(i)} (f_1 (\ldots  f_{n-1}(\tilde {\bf 0})\ldots)) = \tilde e _i  ({\bf{1}}- f_0 (f_1 (\ldots  f_{n-1}(\tilde {\bf 0})\ldots))),
\end{eqnarray*}
where the letter $i$ presents the unique typed $i$ ancestor, so that
\begin{eqnarray*}
	\mathbb{P}(Z_n(i, \cdot) \ne \tilde {\bf 0})= \mathbb E \Big[\mathbb{P}(Z_n(i, \cdot) \ne \tilde {\bf 0}\Big\slash f_0 , f_1,\ldots,f_{n-1})\Big]=\mathbb E [q_{n, i}^{\bf f} ]  .
\end{eqnarray*}
More generally, by the branching property, for any $\tilde z=(z_1, \ldots, z_p) \in \mathbb N^p\setminus \{\tilde {\bf 0}\}$, 
\begin{align}
q^{\bf f}_{n,\tilde z} &:=    \P (|Z_n(\tilde{z},   \cdot)| > 0 \ \slash f_0, \ldots, f_{n-1})  
 =  1 - \prod_{i=1}^p [ 1 - q_{n, i}^{\bf f} ]^{z_i}   
\end{align}
and 
$  
\mathbb{P}(Z_n(\tilde z , \cdot) \ne \tilde {\bf 0})=  	\mathbb E [q_{n, \tilde z}^{\bf f}].$

As in the classical single-type case, the asymptotic behaviour of the quantity $\mathbb E [q_{n, \tilde z}^{\bf f} ]$ above is controlled by the mean matrices and the Hessian matrices of the offspring distributions (see section \ref{Ontheextinctionof}).   

By \cite{FK},  if $
\mathbb E\big[ \ln ^+\vert  M_0 \vert \big] <+\infty,$  then the sequence $\displaystyle  \left(\frac{1}{n} \ln  \vert M_{0, n} \vert \right) _{n\geq 1} $ converges $\mathbb P$-almost surely to some constant limit $\displaystyle \pi_\mu := \lim_{n \to +\infty }   \frac{1}{n} \mathbb{E} \big[ \ln  \vert M_{0, n}\vert \big] $. 
 On the product space ${\mathbb{X}} \times \S^+$, we define the function $\rho $ by setting $\rho (\tilde x, M): = \ln  \vert\widetilde xM\vert$ for $(\tilde x, M) \in {\mathbb{X}} \times \S^+$. This function satisfies the cocycle property, namely for any $M, N \in \S^+$ and $\tilde  x \in  \mathbb{X}$,
\begin{equation} \label{cocycle}
	\rho (\tilde x,MN) = \rho (\widetilde x \, \cdot \, M,  N) + \rho (\tilde x, M).
\end{equation}

Under  hypothesis H3$(\delta$)    introduced below, there exists a unique $\mu$-invariant measure $\nu$ on $\mathbb{X}$ such that for any continuous function $\varphi $ on $\mathbb{X}$,
\[(\mu  * \nu)(\varphi) = \int_{\S^+} {\int_{\mathbb{X}} {\varphi (\tilde x \cdot M)\nu (d\tilde x)\mu ({\rm d}M)}  = }  {\int_{\mathbb{X}} {\varphi (\tilde x)\nu (d\tilde x)}  = } \,\,\nu (\varphi).\]
Moreover, the upper Lyapunov exponent $\pi_\mu$ defined above coincides with the quantity $ \int_{\mathbb{X} \times  \S^+} { \rho (\tilde x, M) \mu ({\rm d}M)\nu (d\tilde  x)} $ and is finite  \cite{BL}.

\vspace{2mm}  
For any $0<\delta <1$, we consider the following    hypotheses concerning   the distribution $\mu$ of the mean matrices  $  M_{n}$  and  the  distributions of the random variables  $\xi_n$  at each step.

\vspace{2mm}

\noindent {\bf Hypotheses    } 

\vspace{2mm}  

H1($\delta$). {\it $\mathbb E[ \vert \ln \mathfrak n(M_1)\vert ^{2+\delta}]< +\infty$.  }

H2. {\it (Strong irreducibility)  There exists no affine subspaces $\mathcal A$ of $\mathbb R ^d$ such that $\mathcal A \cap\mathbb R^p_+$ is non-empty, bounded and invariant under the action of all elements of the support of $\mu$. }

H3($\delta$). {\it The support of $\mu$ is included in $\S^+(B)$ with $B= {1\over \delta}$.}

H4. {\it The upper Lyapunov exponent $\pi_\mu$       equals 0.    }

H5($\delta$). {\it    $ \mu(E_\delta)>0$, where}
$ 
E_\delta:=  \{M \in \S^+\,\, \slash \,\forall \tilde  x \in  \mathbb X, \ \  \ln  \vert\tilde x M\vert\, \ge \delta \}.
 $ 
    \vspace{2mm}

H6. {\it  $\displaystyle \mathbb E\Bigl[{\mu _1\over \vert M_1\vert ^2} (1+\ln^+\vert M_1\vert ) \Bigr] <+\infty$. }

   \vspace{5mm}

 Notice that the moment  hypotheses H1($\delta$),  H3($\delta$) and H6  are satisfied when the offspring generating functions  $f_n$, for $ n \geq 0$,  belong to some $\mathcal G_{\varepsilon, K}$; indeed, in this case,   for any $ 1\leq i, j\leq p,$
 \[ 
2p\varepsilon \leq \vert M_1\vert \leq \sum_{i, j=1}^p\mathbb E [\xi^2_1(i, j) ] \leq p^2K \quad {\rm and } \quad \mu_1\leq p^3 K \qquad  \mathbb P{\rm -a.s.}
 \]

 A lot of researchers investigated the behaviour of the survival probability of  $(Z_n)_{n \geq 0}$  in random environment, under various sets of rather restrictive assumptions.  Following \cite{LPP} and  \cite{DV},    when hypotheses  H1--H6   hold,   for any $\tilde z \in  \mathbb N^p\setminus \{\tilde {\bf 0}\}$ there exist    $ \beta_{\tilde z}>0  $ such that, 
	\begin{equation} \label{extinctionAOP}
	\lim_{n \to +\infty } \sqrt{n} \mathbb{P}(   Z_n(\tilde z, \cdot) \ne \tilde {\bf 0} ) =  \beta_{\tilde z}. 
		\end{equation}
 Notice that  hypothesis H6 above is weaker than the one in \cite{DV}; indeed, the key argument in   \cite{LPP} and  \cite{DV} is based on our  Lemma \ref{serieswidetilde}, which holds under  assumptions H1--H6. 

 The convergence (\ref{extinctionAOP}) relies on a deep understanding, developed in \cite{pham2018},  of the behavior of the 
semi-markovian random walk $(S_n( \tilde x, a))_{n \geq 0}$ defined by $S_n(\tilde x, a)=a+\ln \vert \tilde xM_{0, n}\vert,$ for any $\tilde x \in \X, a \geq 0$ and $n\geq 0$. It is well known  that this Markov walk  satisfies a strong law of large number and a central limit theorem;  denote by   $\sigma ^2:= \displaystyle \lim_{n \to +\infty}  \frac{1}{n} \mathbb E [ S_n^2(\tilde x, a) ] $ its variance  and recall that, under   Hypotheses H1 to H5, the quantity  $\sigma^2$ is positive.

 Here comes the main result of the present paper; it concerns the asymptotic distribution of the random variables $\ln  \vert Z_n (\tilde{z}, \cdot)\vert$ conditioned to non-extinction and requires the strong   assumption that the  offspring distributions $f_n, n \geq 0$,   do belong to some  $\mathcal G_{  \varepsilon, K}$.

 \begin{theo} \label{MAINTHEO}
 Assume that 

(1) there exist $\varepsilon\in ]0, 1[$ and  $K>0$ such that  $\mathbb P$-a.s, the  offspring distributions $f_n,  n \geq 0$, belong to $\mathcal G_{  \varepsilon, K}$;  

(2) there exists $\delta>0$ such that hypotheses {\rm H2}, {\rm H4} and  {\rm H5}($\delta$)  hold.

 Then for any $  \tilde z \in \N^p\setminus\{\tilde{\bf 0} \}$ and $t \ge 0$,
	\begin{equation*}
	\lim_{n \to +\infty }  \P \left( \frac{\ln  \vert Z_n (\tilde{z}, \cdot)\vert}{\sqrt n} \le t \ \Big\slash \ \vert  Z_n (\tilde{z},   \cdot)\vert > 0\right) = \frac{2}{\sigma \sqrt{2 \pi}} \Phi^+ \left( \frac{t}{\sigma} \right),
	\end{equation*}
	where $\Phi^+$ denotes the cumulative function of the Rayleigh distribution:  
	\[\Phi^+ \left(\frac{t}{\sigma} \right) := \int_0^{t } s\exp \left(- \frac{s^2}{2}\right) {\rm d}s.\]  
\end{theo}

The first step to prove this  main theorem is to provide   a limit theorem for the processes $(Z_n(\tilde{z},  j))_{n \geq 0}$, where  $ \tilde z \in \N^p\setminus\{\tilde{\bf 0} \}$ and  $1 \le j \le p$,  conditioned to non-extinction and  randomly  rescaled; this statement is of intrinsic interest and  holds  under weaker the assumptions { \rm H1--H6} (for some $\delta>0$).
\begin{theo}\label{theo2}
Assume that hypotheses { \rm H1--H6} hold for some $\delta>0$. Then, for any $  \tilde z \in \N^p\setminus\{\tilde{\bf 0} \}$ and  $1 \le j \le p$, there exists a probability measure $\nu_{  \tilde z,   j}$ on $\R ^+$ such that for any non-negative continuity point $t$   of the distribution function $s\mapsto \nu_{\tilde z,     j}([0,s])$, 
	\begin{equation*}  
	\lim_{n \to +\infty }  \P \left( \frac{Z_n(\tilde{z},  j)}{\vert  M_{0, n} e_j\vert} \le t \Big\slash \ \vert  Z_n (\tilde{z},   \cdot)\vert > 0 \right) = \nu_{\tilde z,   j} ([0, t]).
	\end{equation*}
	 Furthermore,   if there exist  $\varepsilon\in ]0, 1[$ and $K>0$  such that  $f_n \in \mathcal G_{\varepsilon, K}$ for any $n \geq 0$,   the probability measures $\nu_{\tilde z,   j}$   are supported on $]0, +\infty[$.
\end{theo} 
Theorem \ref{MAINTHEO} is not a direct consequence of Theorem \ref{theo2}; we need an intermediate stage  which  concerns the behaviour of the processes    $(\tilde xM_0\ldots M_{n-1})_{n \geq 0}$, for $ \tilde x \in \X$, conditioned to the event $(\vert  Z_n (\tilde{z},   \cdot)\vert > 0)$. A  close conditioned limit theorem involving the Rayleigh distribution also holds  (see Corollary \ref{theo2pham2018} below) but its  condition is not  the one required here. The following proposition fills this gap and is essential to connect the two statements above.

\begin{prop} \label{proposition2.3}
 Assume that hypotheses { \rm H1--H6}   hold for some $\delta>0$. Then for any $\tilde z \in \N^p\setminus\{\tilde{\bf 0} \}$,  $\tilde x \in \X, a \ge 0$ and  $t \ge 0$,  
	\begin{equation*}
	\lim_{n \to +\infty } \P \left( \frac{S_n (\tilde x, a)}{\sqrt n} \le t \ \Big\slash \ \vert  Z_n (\tilde{z},   \cdot)\vert > 0 \right) = \frac{2}{\sigma \sqrt{2 \pi}} \Phi^+ \left( \frac{t}{\sigma} \right).
	\end{equation*}
\end{prop}
 
 The article is structured as follows. In section  \ref{Auxiliaryresults}, we present some useful auxiliary results  on product of random matrices and properties on varying environment.   Section 
 \ref{Ontherandomenvironment} is devoted to the random environment while the proofs of Theorem \ref{theo2}, Proposition \ref{proposition2.3} and 
Theorem \ref{MAINTHEO} are detailed in sections
\ref{proofoftheorem2}, 
\ref{ProofofProposition2.3}
and 
\ref{proofofmaintheo} respectively.

\vspace{2mm}

 \noindent {\bf Notations.}  {\it  
Let $ (u_n)_{n \geq 0}$ and $ (v_n)_{n \geq 0}$ be two sequences of positive reals; we write

$\bullet\quad $  
	  $u_n \stackrel{c}{\preceq }v_n $ if $u_n \leq c v_n$ for any $n \geq 0$,
	  
	  (and simply $u_n \preceq v_n$ when  $u_n \stackrel{c}{\preceq }v_n $  for some constant $c>0$); 
	  
$\bullet\quad$  $u_n \stackrel{c}{\asymp}  v_n $  when ${1\over c} u_n \leq v_n \leq  c u_n $  for any $n \geq 0$,
	  
	    (and simply $u_n \asymp v_n$ when  $u_n \stackrel{c}{\asymp }v_n $  for some constant $c>0$); 

$\bullet\quad $ $u_n \sim v_n$ if $\displaystyle \lim_{n\to+ \infty} \frac{u_n}{v_n} = 1$.  
}

 \vspace{2mm}

 \noindent {\bf Acknowledegments.} The authors thank  V.A. Vatutin for helpful  comments on the first version of the paper.

\section{Auxiliary results}\label{Auxiliaryresults}
In this section, we state some well known and useful  results  about fluctuations of  products of random matrices  with non-negative entries and some convergence theorems for multi-type Galton-Watson processes in  varying environment. 

\subsection{ On  positive  matrices and their products}

  Following \cite{Hennion}, we endow $\mathbb X$ with a bounded distance ${\mathfrak d}$ such that  any $A\in \mathcal S$ acts on $\mathbb X$ as a contraction with respect to ${\mathfrak d}$.   In the following lemma, we just recall some fundamental  properties of this distance.
\begin{lem} \label{propH}
	There exists a distance ${\mathfrak d}$ on ${\mathbb{X}}$ which is compatible with the standard topology of $\mathbb{X}$ and satisfies the following properties: 
	\begin{enumerate}
		\item $\sup \{  {\mathfrak d}(x,y ) \ \slash \ \tilde x, \tilde y \in  {\mathbb{X}} \} =1 $.
		\item $\vert  x-y\vert \leq 2{\mathfrak d}(x,y)$ for any $\tilde x, \tilde y \in \mathbb X$. 
		\item For any $M \in \mathcal S^+$, set $[ M] := \sup\{{\mathfrak d}(M\cdot x, M\cdot y) \ \slash \ \tilde x, \tilde y \in \mathbb X\}$. Then, 
		\begin{enumerate}
			\item $ {\mathfrak d}(M\cdot x, M\cdot y) \leq [M] {\mathfrak d}(x, y)$ for any $\tilde x, \tilde y \in \mathbb X$;
			
			\item $[MN]\leq [M] [N]$ for any $M, N \in \mathcal S^+$.
		\end{enumerate}
		\item For any $B\geq 1$, there exists $ \rho_B \in ]0, 1[$ such that   $[M]\leq \rho_\delta$ for any $M \in \mathcal  S^+(B)$.
	\end{enumerate}		
\end{lem} 
 Similar statements also hold for the right action of $\mathcal S^+$ and $\mathcal S^+(B)$ on $\mathbb X$.

 The following Property is a direct consequence of Lemma \ref{propH}: up to some normalization, products of matrices in $S^+(B)$ converge to some rank-one matrix. 
 
 For any $M=(M(i, j))_{1\leq i, j\leq p} $ in $S^+ $, we denote by $\overline{M}$ the matrix   with   entries 
 \[
  \overline{M}(i, j)={M(i, j)\over \vert M(\cdot, j)\vert}= {M(i, j)\over \sum_{\ell=1}^k M(\ell, j) }.
  \]

 \begin{property}\label{rank1} 
 	Let  $ {\bf  M}=(M_n)_{n \geq 0} $ be  a sequence of  matrices in $\mathcal S^+(B)$. 
	
	  Then,   the sequence
	$([M_{0, n}])_{n \geq 0}$ converges  exponentially to $0$. In particular, the sequence  $(\overline{M}_{0, n})_{n\geq 0}$  converges  as $n \to +\infty$ towards a rank-one matrix  whose columns are all equal  to ${\bf M}_\infty=({\bf M}_\infty(i))_{1\leq i \leq p}$, where, for any $1\leq j\leq p$, 
	$$
	{\bf M}_\infty(i) :=\lim_{n \to +\infty} {\overline{M}_{0, n}(i, j) }.
	$$ 	
\end{property}
Let us also  recall  some important properties  of matrices  in $\mathcal S^+(B)$.
\begin{lem}\label{keylem}  \cite{FK} Let $T(B)$ be the closed semi-group generated by $\mathcal S^+(B)$. 
	For any $M, N \in T(B)$ and $1 \le i,j,k,l \le p$,
	\begin{eqnarray*} 
	M(i,j) \stackrel{B^2}{\asymp} M(k,l).
	\end{eqnarray*}
	In particular, there exists $ c >1$ such that for any $M, N \in T(B)$ and for any  $\tilde x, \tilde y \in \mathbb{X}$,
	\begin{enumerate}
		\item  $ \left\vert { M x} \right\vert \stackrel{c}{\asymp} \left\vert M \right\vert \,\, \mbox{and} \,\, \vert\tilde y M\vert \stackrel{c}{\asymp} \vert M\vert$, 
		\item $\vert \widetilde yMx\vert \,\, \stackrel{c}{\asymp} \,\,\vert M\vert $,
		\item $\vert MN\vert  \stackrel{c}{\asymp} \vert M\vert \vert N\vert $.
	\end{enumerate}

\end{lem} 
\subsection{ Limit theorem for products of random positive matrices}

Throughout this  subsection, the matrices $M_n, n \ge 0,$  are i.i.d.   and their law $\mu$ satisfies hypotheses H1--H5  for some $\delta >0$. We introduce the homogenous Markov chain $(X_n)_{n \geq 0}$ on $\mathbb{X}$ defined by the initial value $X_0=\tilde x  \in \mathbb{X}$ and for $n \ge 1$,
\[ X_{n} =  \tilde x \cdot M_{0, n}.\]
Its transition  probability $P$  is given  by: for any $\tilde  x \in  \mathbb{X}$ and any bounded Borel function $\varphi:  \mathbb{X} \to \mathbb{R}$,
\[P\varphi (\tilde x) := \int_{\S^+} {\varphi (\tilde x \cdot M) \mu ({\rm d}M)}. \]

In the sequel, we are interested in the left linear action $ \tilde x\mapsto \widetilde x{M_{0, n}}$ of the right products $M_{0, n}$, for any $\tilde  x \in  \mathbb X$.   By simple transformation, we see that
\[\widetilde x{M_{0, n}}= e^{\ln  \vert \widetilde x{M_{0, n}}\vert}\   \widetilde x\cdot {M_{0, n}},\] 
which turns it natural to consider the random process $(S_n)_{n \geq 0}$   defined by: for any $\tilde  x \in  \mathbb X, a \in \mathbb R$ and  $n \geq 1$, 
\[S_0=S_0(\tilde x, a) :=a\quad {\rm and} \quad  {S_n} =S_n(\tilde x, a):= a+\ln  \vert \widetilde x{M_{0, n}}\vert .
\]
In order to simplify the notations,    let $S_n(\tilde x):= S_n(\tilde x, 0)$ for any $\tilde x\in \mathbb X$ and any $n \geq 0$. 
By iterating the cocycle property (\ref{cocycle}), the basic decomposition of $S_n(\tilde x, a)$ arrives: 
\begin{equation*}
	S_n (\tilde x, a) =a+ \ln  \vert \tilde x {M_{0, n}}\vert  = a+ \sum_{k = 0}^{n-1} {\rho ({X_k}, { M _{k}})}.
\end{equation*}
 It is noticeable that for any $a \in \mathbb{R}$, the sequence $(X_n, S_n)_{n \geq 0}$ is a Markov chain on $\mathbb{X} \times \mathbb{R}$ whose transition  probability $\widetilde P$ is defined by: for any $(\tilde x, a) \in \mathbb{X} \times \mathbb{R}$ and any bounded Borel function $\psi: \mathbb{X} \times \mathbb{R} \to \mathbb{C},$
\[\widetilde P\psi (\tilde x, a) = \int_{\S^+} {\psi (\tilde x \cdot M,a+\rho(\tilde x, M)) \mu ({\rm d}M)}. \]
For any $\tilde  x \in  \mathbb{X}$ and $a\geq0$, we denote by $\mathbb{P}_{\tilde x, a} $ the  probability measure on $(\Omega, \mathcal{F}, \mathbb{P})$ conditioned to the event $(X_0 = \tilde x, S_0 =a)$ and  $\mathbb{E}_{\tilde x, a} $ the corresponding expectation;   the index $a$ is omitted when $a =0$  and $\mathbb P  _{\tilde x}$ denotes  the corresponding  probability.

We set $\mathbb R_*^+:= \mathbb R^+{\setminus}\{0\}$ and  define  $\widetilde P_+$ the restriction of $\widetilde P$ to $\mathbb{X} \times \mathbb{R}^+$:  for $a>A$ and any $\tilde  x \in  \mathbb{X},$ 
\[\widetilde P_+ ((\tilde x, a),\cdot) = {\bf 1}_{\mathbb{X} \times \mathbb{R}^+}(\cdot)\widetilde P((\tilde x, a),\cdot). \]
 Furthermore, we introduce  

$\bullet$  the first (random) time  at which  the random process $(S_n(\tilde x, a))_{n \geq 0}$ becomes non-positive:
\begin{eqnarray*}  
\tau = \tau_{\tilde x, a}:= \min \{ n \ge 1\ \slash \, S_n(\tilde x, a)  \le 0\}. 
\end{eqnarray*}

$\bullet$ 
the  minimum  $m_n$ for $n \geq 1$, defined by
\[
m_n = m_n(\tilde x, a) := \min\{S_1(\tilde x, a),\ldots,S_n(\tilde x, a) \},
\]
 and we set  
\[{\bf m}_n(\tilde x,  a): = \mathbb{P}_{\tilde x, a}({m_n} > 0)= \mathbb P_{\tilde x, a}(\tau>n)=\mathbb P(\tau_{\tilde x, a}>n).\] 
In order to simplify the notations,  when $a=0$   let $\tau_{\tilde x} := \tau_{\tilde x, 0}, \ m_n(\tilde x) :=m_n(\tilde x, 0)$ and  ${\bf m}_n(\tilde x) :={\bf m}_n(\tilde x,  0)$ for any $\tilde x\in \mathbb X$ and any $n \geq 0$.

 We recall some important results about the behaviour of the random products of variables $M_{0, n}$ and the distribution of $\tau$, under $
\mathbb P_{\tilde x, a}$.

 \begin{prop} \label{prop1pham}  \cite{pham2018}
Assume hypotheses {\rm H1 -- H5}  for some $\delta >0$.	Then for any $\tilde x \in  \mathbb X$ and $a \ge 0$, the sequence $\Bigl(\mathbb E _{\tilde x, a}[ S_n; \tau  >n ] \Bigr) _{n \ge 0}$ converges to some quantity $	V(\tilde x, a)$.
The function $V$ is $\widetilde P_+$-harmonic on $\mathbb X\times \mathbb R^+$ and satisfies the following properties: 

	\begin{enumerate}
	 		
		\item for  any $\tilde x \in  \mathbb X$, the function $V(\tilde x, \cdot)$ is increasing on $\mathbb R^+ $;
				
		\item there exist $c >0$ and $A >0$ such that for any $\tilde x \in  \mathbb  X$ and $a\ge 0$,  
		\[
		\frac{1}{c} \vee (a-A) \leq V(\tilde x, a)  \leq c(1+a);
		\]
		
		\item for any $\tilde x \in  \mathbb X$, the function $V(\tilde x, .)$ satisfies $\displaystyle \lim_{a \to +\infty} \frac{V(\tilde x, a) }{a} =1$.

	\end{enumerate}
\end{prop}
The next statement allows to control the tail of the distribution of the random variable  $\tau$.
\begin{theo} \label{theo1pham}   \cite{pham2018}
	Assume  hypotheses { \rm H1--H5}  for some $\delta >0$. Then for any $\tilde x \in  \mathbb X$ and $a \ge0$,
	\[
	\mathbb P _{\tilde x, a}(\tau >n) \sim \frac{2 }{\sigma \sqrt{2 \pi n}} V(\tilde x, a)\,\,\mbox{as} \,\, n \to +\infty,
	\]
	where $\sigma ^2>0$ is the variance of the Markov walk $(X_n, S_n)_{n \geq 0}$.
	Moreover, there exists a constant $c$ such that for any $\tilde x \in  \mathbb X$, $a\ge 0$ and $n \ge 1$,
	\[
	0\leq \sqrt n \mathbb P _{\tilde x, a}(\tau  > n) \leq c V(\tilde x, a).
	\]	
\end{theo}

  Our hypotheses { \rm H1--H5}  correspond  to  hypotheses{ \rm P1--P5}    in \cite{pham2018}, except  that hypothesis H1 is weaker than P1. Indeed, the existence of moments of order $2+\delta$  suffices.  This ensures    that the   map $t \mapsto P_t$  in Proposition 2.3 of \cite{pham2018}  is  $C^2$,   which suffices for this Proposition  to hold. Moreover, the martingale  $(M_n)_{n \geq 0} $ which approximates the process  $(S_n(x))_{n \geq 0}$ belongs to $\mathbb L^p$ for $p=2+\delta$ (and not for any $p>2$ as stated in \cite{pham2018} Proposition 2.6). This last property was useful in \cite{pham2018} to achieve the proof of Lemma 4.5, by choosing $p$ great enough in such a way that $(p-1)\delta -{1\over 2}>2\varepsilon$ for some fixed constant $\varepsilon>0$. Recently, by following the same strategy as C. Pham, M. Peign\'e and W. Woess    improved this part of the proof,  by allowing various  parameters (see \cite{PW},  Proof of Theorem 1.6 (d)).

As a direct consequence,   up to some normalisation and conditioned to the event $[m_n>0]$,  the sequence $\left( S_n \right)_{n \ge 1}$ converges weakly to the Rayleigh  distribution.  
\begin{coro} \label{theo2pham2018}  \cite{pham2018}
	Assume that hypotheses  { \rm H1--H5} hold for some $\delta>0$. Then  for any $\tilde x \in  \mathbb X$, $a\ge 0$ and $t >0$, 
	\[
	\lim_{n \to +\infty} \mathbb P _{\tilde x, a}\left( \frac{S_n}{ \sqrt n} \leq t\Big\slash \tau >n\right)  =
	\Phi^+\left({t\over \sigma}\right).
		\]

\end{coro}

Theorem \ref{theo1pham}  leads to some upper bound in the local limit theorem. 

\begin{coro}\label{theolocal}
Assume hypotheses  { \rm H1--H5} for some $\delta>0$. Then there exists a positive constant  $c$ such that for any $\tilde x \in  \mathbb X, a, b \ge0$ and $n \geq 1$,
\[
	0\leq \mathbb P _{\tilde x, a} \big(S_n \in [b, b+1[, \tau  > n \big) \leq c {(1+a)(1+b)\over n^{3/2} }.
\]
\end{coro}
{\bf Proof.}  We follow and adapt   the strategy of  Proposition 2.3 in \cite{ABKV}.  For fixed $\tilde x \in \mathbb X $ and $a, b>0$, we write 
\begin{align*}
E_n=E_n(x, a, b)&:=\Bigl(S_n(\tilde x, a) \in [b, b+1[, \ \tau_{\tilde x, a}  > n\Bigr)\\
&=
\Bigl(e^a\vert \tilde x M_{0, n}\vert \in   [e^b, e^{b+1}[, \ e^a\vert \tilde x M_0\vert > 1, \ldots, \ e^a \vert \tilde x M_{0, n}\vert > 1\Bigr)
\\
&\subset
  \Bigl( \tau_{\tilde x, a}  > n/3\Bigr) 
\\
 & \qquad   \cap \Bigl(\vert \tilde xM_{0, n}\vert \in   [e^{b-a},  e^{b-a+1}[
 , \ e^a\vert \tilde x M_{0, [2n/3]+1}\vert > 1, \ldots, \ e^a \vert \tilde x M_{0, n}\vert > 1
 \Bigr). 
\end{align*}
Let us decompose $M_{0, n}$ into three parts,  using  the notation $M_{k, n}= M_{k}\ldots M_{n-1} $ for any $0\leq k < n$ and $n \geq 1$. 
It holds that $M_{0,n}=M_n'M_n''M_n'''$ with 
\[
M_n'= M_{0,[n/3]}= M_0 \ldots M_{[n/3]-1},\quad   M_n''= M_{[n/3], [2n/3]}=M_{[n/3]}\ldots M_{[2n/3]-1}\]
and 
\[ M_n'''=M_{[2n/3], n}=
M_{[2n/3]}\ldots M_{n-1}.
\] 
By  Lemma  \ref{keylem}, we may write, on the one hand 
\begin{align*}
\Bigl(\vert \tilde xM_{0, n}\vert \in   [e^{b-a},  e^{b-a+1}[\Bigr)
&
=\Bigl(\vert \tilde xM_n'M_n''M_n'''\vert \in   [e^{b-a},  e^{b-a+1}[\Bigr)
\\
&
{\stackrel{\mathbb P \text{-a.s.}}{ \subset}}
\Bigl(\vert  M_n'M_n''M_n'''\vert \in   [e^{b-a},  ce^{b-a+1}[\Bigr)
\\
&
{\stackrel{\mathbb P \text{-a.s.}}{ \subset}}
\Bigl( \vert M_n''\vert   \in   \Bigl[{e^{b-a}\over \vert  M_n'\vert \ \vert  M_n'''\vert},  c^3{e^{b-a+1}\over \vert  M_n'\vert \ \vert  M_n'''\vert}\Bigl[\Bigr)
\\
&
{\stackrel{\mathbb P \text{-a.s.}}{ \subset}}
\Bigl( \vert\tilde x M_n''\vert   \in   \Bigl[{e^{b-a}\over c\vert  M_n'\vert \ \vert  M_n'''\vert},  c^3{e^{b-a+1}\over \vert  M_n'\vert \ \vert  M_n'''\vert}\Bigl[\Bigr)
\end{align*}
and on the other hand,      for any $[2n/3]<k\leq n$,
\[{1\over c} \vert \tilde xM_{0, k}\vert\  \vert M_{k+1, n}\vert\leq \vert  \tilde x M_{0, n}\vert
\ \mathbb P\text{-a.s.  }
\]
This yields that
\begin{align*}
&\Bigl(\vert \tilde xM_{0, n}\vert \in   [e^{b-a},  e^{b-a+1}[
 , \  e^a\vert \tilde x M_{0, [2n/3]+1}\vert > 1, \ldots, e^a \vert \tilde x M_{0, n}\vert > 1
 \Bigr)\\
 &\qquad \qquad \qquad \qquad \qquad {\stackrel{\mathbb P \text{-a.s.}}{ \subset}}\quad \Bigl(\vert M_{[2n/3]+1, n}\vert <ce^{b+1},   \ldots,  \vert  M_{n-1, n}\vert <ce^{b+1}\Bigr)\\
&\qquad \qquad \qquad\qquad \qquad {\stackrel{\mathbb P \text{-a.s.}}{ \subset}}\quad \Bigl(\vert \tilde x M_{[2n/3]+1, n}\vert <ce^{b+1},   \ldots,  \vert \tilde x M_{n-1, n}\vert <ce^{b+1}\Bigr).
\end{align*}
Finally,
\[
E_n \quad  {\stackrel{\mathbb P \text{-a.s.}}{ \subset}}\quad E'_n \cap E''_n \cap E'''_n 
\]
with
\[E'_n=E'_n(x, a, b)= \Bigl( \tau_{\tilde x, a}  > n/3\Bigr), 
  \quad  E''_n=E''_n(x, a, b)= \Bigl( \vert\tilde x M_n''\vert   \in   \Bigl[{e^{b-a}\over c \vert  M_n'\vert \ \vert  M_n'''\vert},  c^3{e^{b-a+1}\over \vert  M_n'\vert \ \vert  M_n'''\vert}\Bigl[ \ \Bigr)\]
and
\[
E'''_n=E'''_n(x, a, b)= \Bigl(\vert \tilde x M_{[2n/3]+1, n}\vert <ce^{b+1},   \ldots,  \vert \tilde x M_{n-1, n}\vert <ce^{b+1}\Bigr).\]
The events  $E'_n$ and $E'''_n$ are  measurable with respect to the $\sigma$-field  $\mathcal T_n$ generated by $M_0, \ldots, M_{[n/3]-1}$ and $M_{[2n/3] }, \ldots, M_{n-1}$; consequently,
\begin{align*}
\mathbb P(E_n)&\leq \mathbb P\Bigl(\mathbb P \big(E'_n\cap E''_n\cap E'''_n \ \vert \ \mathcal T_n \big)\Bigr)
\\
&=\mathbb E \big[ {\bf 1}_{E'_n\cap E'''_n  } \mathbb P( E''_n \ \vert \ \mathcal T_n ) \big].
\end{align*}
The random variable $\ln  \vert \tilde x M''_n\vert$ are independent on  $\mathcal T_n$ and their distribution coincides with the one of $ S_{[n/3]}(\tilde x, 0)$. Therefore, by the classical local limit theorem for product of random matrices  with non-negative entries,  for any  $\tilde x \in \X$ and $a, b \ge 0$, 
\begin{align*}
\mathbb P \big( E''_n \ \vert \ \mathcal T_n \big)&
{\stackrel{\mathbb P \text{-a.s.}}{ \leq}} 
\sup_{A\in \mathbb R} \mathbb P\Bigl( \vert\tilde x M_n''\vert   \in   \Bigl[A,  c^4 e A\Bigl[\Bigr) \preceq {1\over \sqrt{n}}.
\end{align*}
 Since the events  $E'_n$ and $E'''_n$ are independent,  it follows that 
 \begin{equation}\label{probEn}
 \mathbb P(E_n)\preceq {\mathbb P(E'_n\cap E_n'')\over \sqrt{n}}= {\mathbb P(E'_n)\mathbb P(E_n'')\over \sqrt{n}}.
 \end{equation}
The probability of $E'_n(x, a, b)$ is controlled by Theorem \ref{theo1pham}:  uniformly   in $\tilde x \in \X$ and $ a, b \ge 0$, 
 \begin{equation}\label{probE'n}
 \mathbb P(E'_n(x, a, b)) \preceq {1+a\over \sqrt{n}}.
 \end{equation} 
To control the probability of the event $E_n'''(x, a, b) $, we introduce the stopping time  
\[ \tau^+ := \min \{ n \ge 1 \ \vert \, S_n  \geq 0\} 
\]
and notice that its  distribution satisfies the same tail condition as $\tau$.  Since  
\[
(M_{[2n/3] }, M_{[2n/3]+1}, \ldots,  M_{n-1} )  
{\stackrel{  \text{dist}}{=}}
(M_ 0, M_1 , \ldots, M_{n-[2n/3]-1}), \]  
it holds uniformly in $\tilde x, a, b$ that, 
\begin{align}\label{probE'''n}
\mathbb P(E_n'''(x, a, b)) &\leq \mathbb P\Bigl(\vert \tilde x M_{0}\vert <ce^{b+1},   \ldots,  \vert \tilde x M_{0,n- [2n/3]}\vert <ce^{b+1}\Bigr)
\notag \\
&=
\mathbb P_{\tilde x, -\ln c-b-1}(\tau^+>n-[2n/3]) \notag\\
&\preceq  {1+b\over \sqrt{n}}. \end{align}
The proof is done, by combining (\ref{probEn}), (\ref{probE'n}) and (\ref{probE'''n}).

\rightline{$\Box$}

\subsection{On the probability of extinction in varying environment} \label{Ontheprobabilityofextinctioninvaryingenvironment}

	In this section we state some useful  results concerning  multi-type Galton-Watson processes in varying environment  ${\bf g}= (g_n)_{n \geq 0}$.

		For any $1\leq j \leq p$,   the quantity 
		 $\left\vert M^{{\bf g}}_{0, n}(\cdot, j) \right \vert$ equals the  mean number  
				 $\E[  Z_n(\tilde {\bf 1}, j)]$ of   particles of type $j$ at generation $n$,  given that there is one ancestor of each type at time $0$.		
	By  Lemma \ref{keylem},  if  all the $M_{g_n}$  belong to $ \S^+(B)$, then    $\displaystyle  \left\vert M^{{\bf g}}_{0, n}(\cdot, j) \right \vert\asymp  \left\vert M^{{\bf g}}_{0, n}(i, j)  \right \vert\asymp  \left\vert M^{{\bf g}}_{0, n}  \right \vert$ for any  $1\leq i, j\leq p$ and $n \geq 0$.  Furthermore,  by Property  \ref{rank1}, the sequence  of normalized matrices $(\overline{M}^{{\bf g}}_{0, n})_{n \geq 0}$ converges as $ n \to +\infty$ towards a rank-one matrix   with common column vectors  ${\bf M}_\infty^{\bf g}=({\bf M}_\infty^{\bf g}(i))_{1\leq i \leq p}$ .
	
		The following statement brings together several results  obtained by D. Dolgopyat and al. \cite{DHKP},  O. D. Jones \cite{Jones}  and G. Kersting  \cite{kersting2017} in the varying environment framework. 
		The last point of this statement is a new key to the main   theorem of this paper: the conditioned central limit theorem  for multi-type Galton-Watson processes in random environment.

		\begin{prop}\label{varyingenvi} Let $Z^{\bf g}=(Z^{\bf g}_n)_{n \geq 0}$ be a  $p$ multi-type Galton-Watson process in varying environment ${\bf g} = (g_n)_{n \geq 0}$. 
		
		$\bullet$ Assume that  there exists $B>1$  such that for any $n \geq 0$,  the mean matrices $M_{g_n}$ belong to $\S^+(B)$. 

		Then,
			
			1.  if   for some (hence every)  $i, j \in \{1, \cdots, p\}$,
				\begin{equation} \label{conditionjones}
					\sum_{n \geq0}   {1\over \vert M^{\bf g}_{0, n}\vert }{ \sigma^2_{g_n}(i,j)\over \vert M_{g_{n+1}}\vert ^2}< +\infty,
				\end{equation}
				then there exists a non-negative  random column vector ${\mathcal W}^{\bf g}= (\mathcal W^{\bf g}(i))_{1\leq i \leq p} $ such that $\mathbb E [\mathcal W^{\bf g}(i )] = {\bf M }_\infty^{\bf g}(i)$ and as $n \to +\infty$,  for every $i, j \in \{1, \ldots, p\}$, 
				\begin{equation} \label{w(i)bfg}
				 	\mathcal W^{\bf g}_n (i, j):=	{Z^{\bf g}_n(i, j) \over\left\vert M^{{\bf g}}_{0, n}(\cdot, j) \right \vert } ={Z^{\bf g}_n(i, j) \over \vert \mathbb E[ Z^{\bf g}_n(\tilde {\bf 1}, j)] \vert}\quad \xrightarrow{\mathbb L^2(\P _{{\bf g}}) } \quad  \mathcal W^{\bf g}(i).
				\end{equation} 
				
				$\bullet$  If it is further assumed that there exists $\varepsilon, K>0$ such that all  the  generating functions $g_k, k\geq 0$, belong to $\mathcal G_{\varepsilon, K}$, then, 
				
				2.  the extinction of the process $Z^{\bf g}$  occurs with probability  $ q^{\bf g}_{  i}<1$ for some (hence  every) $i\in \{1, \ldots, p\} $    if and only if  for some (hence every) $ i, j\in \{1, \ldots, p\} $, 
				\begin{equation}\label{meandivergence}
					\sum_{n\geq 0}{1\over M^{{\bf g}}_{0, n}(i, j)}<+\infty.
				\end{equation} 
			
			3. 
				under conditions (\ref{meandivergence}) and (\ref{conditionjones}),  for any  $1\leq i\leq p,  $ it holds
				\begin{equation}\label{extinction-g}
			  \ \mathbb P\bigl(\mathcal W^{\bf g}(i)=0)  = \mathbb P\bigl(\liminf_{n} \vert  Z^{\bf g}_n(i, \cdot)\vert =0\bigr)=q^{\bf g}_i .
				\end{equation}
		 
		\end{prop}
	Let us shortly comment on this statement. 
		
$\bullet$	  For a single-type  supercritical Galton-Watson process $(Z_n)_{n \ge 0}$  in constant environment, it is well known that the sequence 
	$(\mathcal W_n)_{n \geq 0}$, where $\mathcal W_n:= Z_n/\E[Z_n], $ is a non-negative martingale, hence it converges $\mathbb P$-a.s. towards some non-negative limit $\mathcal W$ (and in ${ \mathbb L}^1$ under some other moment conditions). The first assertion corresponds to a weak version of this  property for multi-type Galton-Watson processes in  varying environment, without the martingale's argument which fails here.

	$\bullet$	  The second assertion  means that condition (\ref{meandivergence}) is equivalent to the fact that $Z^{\bf g}$ is supercritical.

$\bullet$  The third assertion corresponds to the famous ``Kesten-Stigum's theorem''; this is new result in the context of multi-type Galton-Watson processes in varying environment, the proof is detailed in subsection \ref{Ontheprobabilityofextinctioninvaryingenvironment}. We refer   to  \cite{kersting2017} and references therein.

\noindent {\bf Proof. }
 Throughout this proof, in order to simplify the notations, we omit the exponent $\bf g$, except  at the end when we need to specify the environment.
 
 (1)  Assertion 1 corresponds to Theorem 1 in \cite{Jones}; the fact that  the mean matrices $M_{g_n}, n \geq 0$, belong to    $\S^+(B)$ readily implies that they are ''allowable and weakly ergodic'' in the sense of  O. D. Jones  \cite{Jones}.

 (2) Assertion 2  follows by combining Proposition  2.1 (e)  and Theorem 2.2 in \cite{DHKP}.    As far as we know,  without the restrictive assumption $g_n \in \mathcal G_{\varepsilon, K}$ for all $n \geq 0$, there exists no    criteria in the literature  in terms of the mean matrices ensuring the super-criticality of the  process $Z^{\bf g}$.

(3) In \cite{Jones}, the author establishes  some conditions on ${\bf g}$ which ensure that  equality  (\ref{extinction-g})  holds; nevertheless, as claimed there, ``it can be difficult to check them'', except in some restrictive cases. Therefore, as far as we know,  Assertion 3 is  a new statement;   we  detail its proof here, by following the strategy  developed by G. Kersting  (Theorem 2 (ii) in \cite{kersting2017}) and by using some estimations obtained in \cite{DHKP}.

 By construction of the $\mathcal W^{\bf g}(i), 1\leq i \leq p$, the inclusion
 \[ (\liminf_{n} \vert  Z_n^{\bf g}(i, \cdot)\vert =0)\subset (\mathcal W^{\bf g}(i)=0)  \]
 is obvious.
Hence, it suffices   to show that  $\mathbb P( \mathcal W^{\bf g}(i)  =0, \liminf_n\vert Z^{\bf g} _n(i, \cdot)\vert \ \geq 1)=0$.   We decompose the argument into two steps.

\noindent 
{\bf  Step 1.}
{\it Comparison between  $\mathbb P(\liminf_n\vert Z^{\bf g}_n(\ell, \cdot)\vert \ =  0)$ and  $ \mathbb P(\mathcal W^{\bf g} (i)=0 )$.}

By formula (7) in \cite{DHKP},  since the sequence ${\bf g}_1=(g_k)_{k \geq 1}$ belongs to $\mathcal E$,  the functions $g_{1, n}=g_1\circ\ldots \circ g_n$ satisfy the following property: for $s   \in [0, 1[^p$,  
 \begin{equation} \label{generating}
\sum_{k=1}^n{1\over \vert M^{\bf g}_{1, k} \vert} \preceq {1\over \vert {\bf  1}-g_{1, n}(s)\vert}
\preceq {1\over  \vert M^{\bf g}_{1, n} \vert} {1\over \vert{\bf 1}-s \vert }  +\sum_{k=1}^n{1\over \vert M^{\bf g}_{1, k} \vert}.
\end{equation}
By  convexity of $g_0$,  there exists $c_0\geq 1$  such that for any $1\leq i\leq p$ and $t   \in [0, 1[^p$, 
\begin{equation}\label{convexg0}
 {1\over c_0}\leq {1-g_0^{(i)}(t)\over \vert 1-t\vert}    \ 
\leq c_0.
\end{equation}(We detail  the argument at the end  of the present proof).
Thus, by combining   (\ref{generating}), (\ref{convexg0}) and Lemma \ref{keylem},   for any $1\leq i\leq p$, it holds that
 \begin{equation*}  
\sum_{k=1}^n{1\over \vert M^{\bf g}_{0, k} \vert} \preceq {1\over 1-g^{(i)}_{0, n}(s) }
\preceq {1\over  \vert M^{\bf g}_{0, n} \vert} {1\over \vert{\bf 1}-s \vert }  +\sum_{k=1}^n{1\over \vert M^{\bf g}_{0, k} \vert}.
\end{equation*}

This yields 

- on the one hand,  by  choosing $s={\bf 0}$,  
\[
 {1\over \displaystyle    \mathbb P(\vert Z^{\bf g}_n(i, \cdot)\vert \geq 1)} 
\succeq  \sum_{k=1}^n{1\over \vert M^{\bf g}_{0, k} \vert};
\]

- on the other hand,      by setting $s_\lambda= (e^{-\lambda/\vert M^{\bf g}_{0, n}(\cdot, 1)\vert},  1, \ldots, 1 )$ with   $\lambda >0$,
\[
 {1\over   1-\mathbb E\left[ e^{-\lambda {Z^{\bf g}_n(i, 1)\over \vert M^{\bf g}_{0, n}(\cdot, 1)\vert}}\right] }= {1\over \displaystyle  1-g_{0, n}^{(i)}(s_\lambda)}\leq 
 {1\over  \vert M^{\bf g}_{0, n} \vert} {1\over1-e^{-{\lambda \over \vert M^{\bf g}_{0, n}(\cdot, 1)\vert}}}   +\sum_{k=1}^n{1\over \vert M^{\bf g}_{0, k} \vert}.
\]

This readily implies that for any $1\leq i, \ell \leq p$,  
\[
{1\over \displaystyle  1-\mathbb E\left[ e^{-\lambda {Z^{\bf g}_n(i, 1)\over \vert M^{\bf g}_{0, n}(\cdot, 1)\vert}}\right] }
\preceq  {1\over  \vert M^{\bf g}_{0, n} \vert} {1\over1-e^{-{\lambda \over \vert M^{\bf g}_{0, n}(\cdot, 1)\vert}}}  
+{1\over \displaystyle   \mathbb P(\vert Z^{\bf g}_n(\ell, \cdot)\vert \geq 1)}. 
\]  
By Lemma \ref{keylem}, it holds  that $\vert M_{0, n}(\cdot, 1)\vert \asymp  \vert M^{\bf g}_{0, n} \vert $; furthermore,   these quantities  tend  to $+\infty$  as $n\to +\infty$,     by (\ref{meandivergence}). Hence
\[ 
{1\over \displaystyle  1-\mathbb E\left[ e^{-\lambda  \mathcal W^{\bf g} (i)}\right]}
\preceq  {1\over   \lambda}  
+{1\over \displaystyle   \mathbb P(\liminf_n\vert Z^{\bf g}_n(\ell, \cdot)\vert \geq 1)}. 
\]
Letting $\lambda\to +\infty$ yields that $\quad \displaystyle {1\over 1-\mathbb P(\mathcal W^{\bf g} (i)=0)}
\preceq
 {1\over \displaystyle   \mathbb P(\liminf_n\vert Z^{\bf g}_n(\ell, \cdot)\vert \geq 1)}. 
$
	
 In other words, there exists a constant $\kappa\geq 1$ such that  for any $i, \ell \in \{1, \ldots, p\}$,
\[
 \mathbb P(\liminf_n\vert Z^{\bf g}_n(\ell, \cdot)\vert \geq 1)\ \leq \ 
 \kappa  \mathbb P(\mathcal W^{\bf g} (i)>0),  
\]
which implies that when $\mathbb P(\mathcal W^{\bf g} (i)>0) \leq {1\over 2\kappa}$,   \begin{equation}\label{probas}
 \mathbb P(\liminf_n\vert Z^{\bf g}_n(\ell, \cdot)\vert \ = 0)\ \geq \ 
\Bigl( \mathbb P(\mathcal W^{\bf g} (i)=0\Bigr)^{2\kappa}.  
\end{equation}
To get  this last inequality, we use the following elementary lemma. 

\begin{lem} \cite{kersting2017}
Let  $\kappa \geq 1$ and $A, B$ be two events such that $\mathbb P(A)\leq \kappa \mathbb P(B)$ and $\mathbb P(B) \leq {1\over 2\kappa}$. Then 
\[\mathbb P(\overline{A})\geq \mathbb P\mathbb (\overline B)^{2\kappa}.\]
\end{lem}

\noindent {\bf  Step 2.} {\it A martingale argument.}

As G. Kersting in \cite{kersting2017},  we  introduce a martingale defined by: for $k\geq 0$ and any $i \in \{1, \ldots, p\}$,  
\[
{\mathcal M} _k:= \mathbb P(\mathcal W^{\bf g} (i)=0\Big\slash Z^{\bf g}_0, \ldots, Z^{\bf g}_k).
\]
It is known that ${\mathcal M }_k\to {\bf 1}_{(\mathcal W^{\bf g} (i)=0)}\ \mathbb P$-a.s.  as $k \to +\infty$ by standard martingale theory; in particular it converges $\mathbb P$-a.s. towards $1$ on the event $(\mathcal W^{\bf g} (i)=0)$.

The branching property of the process $(Z^{\bf g}_n)_{n \geq 0}$ is used to express ${\mathcal M } _k$ in another  form. It is noticeable that $\mathcal W^{\bf g}(i)$ depends on the whole sequence ${\bf g}$;  let us set  $ {\bf g}_k:= (g_l)_{l\geq k}$ and denote $\mathcal W^{{\bf g}_k}(i)$  the random variable defined as in (\ref{w(i)bfg}) but with respect to the Galton-Watson process $  Z^{{\bf g}_k}$  corresponding to the  environment  ${\bf g}_k $.  By the branching property, 
\[
{\mathcal M} _k= \mathbb P(\mathcal W^{{\bf g}_k}(1)=0)^{Z^{\bf g}_k(i, 1)} \times   \ldots\times  \mathbb P(\mathcal W^{{\bf g}_k}(p)=0)^{Z^{\bf g}_k(i, p)} , 
\]
so that  for $1\leq  \ell \leq p,$ as $k \to +\infty$, 
\begin{equation}\label{convergence}
\mathbb P(\mathcal W^{{\bf g}_k}( \ell)=0)^{Z^{\bf g}_k(i,  \ell)} \longrightarrow 1 \quad  \mathbb P\text{-a.s. on the event} \quad  (\mathcal W^{\bf  g}(i)=0). 
\end{equation}
The same property holds, replacing the event $(\mathcal W^{{\bf g}_k}( \ell)=0)$ by $(\liminf_n\vert Z^{{\bf g}_k}_n(\ell, \cdot)\vert \ = 0)$, namely: for any  $\ell \in \{1, \ldots, p\}$, as $k \to +\infty$,
\begin{equation}\label{convergencebis}
 \mathbb P(\liminf_n\vert Z^{{\bf g}_k}_n(\ell, \cdot)\vert \ = 0)^{Z ^{\bf g}_k(i,  \ell)} \longrightarrow 1 \quad  \mathbb P\text{-a.s. on the event} \quad  (\mathcal W^{\bf  g}(i)=0).
\end{equation}
Indeed,  every subsequence $\Bigl(\mathbb P(\displaystyle \liminf_n\vert Z^{{\bf g}_{k_r}}_n(\ell, \cdot)\vert \ = 0)^{Z^{\bf g}_{k_r}(i,  \ell) }\Bigr)_{r\geq 0}$ has a further subsequence which converges to $1$.
In order to apply inequality (\ref{probas}), we distinguish two cases.

(i) Either $\mathbb P(\mathcal W^{{\bf g}_{k_r}}( \ell)>0) \leq {1\over 2\kappa}$  for  $k$ large enough;   we may apply (\ref{probas}) and (\ref{convergence}) to obtain, $\mathbb P$-a.s.  on $ (\mathcal W^{\bf  g}(i)=0)$,
 \begin{equation*} 
\liminf_{r \to +\infty} \mathbb P(\liminf_n\vert Z^{{\bf g}_{k_r}}_n(\ell, \cdot)\vert \ = 0)^{Z^{\bf g} _{k_r}(i,  \ell) }
\geq 
\liminf_{r \to +\infty} \mathbb P(\mathcal W^{{\bf g}_{k_r}}( \ell)=0)^{2\kappa Z ^{\bf g}_{k_r}(i,  \ell) }=1.
\end{equation*}

(ii) Or there exists   a further  subsequence $(k'_r)_{k \geq 0}$  such that  
$\mathbb P(\mathcal W^{{\bf g}_{k'_r}}( \ell)>0)>{1\over 2\kappa}$. Hence, (\ref{convergence}) implies that $Z^{\bf g}_{k'_r}(i, \ell) \to 0\quad \mathbb P$-a.s.  on  $ (\mathcal W^{\bf  g}(i)=0)$,  as $r \to +\infty$; in other words $Z^{\bf g}_{k'_r}(i, \ell) = 0$ for $r$ large enough,     thus $\mathbb P(\displaystyle \liminf_n\vert Z^{{\bf g}_{k'_r}}_n(\ell, \cdot)\vert \ = 0)^{Z^{\bf g}_{k'_r}(i,  \ell) }=1.$

\vspace{2mm}
\noindent Finally, in both cases, convergence  (\ref{convergencebis}) holds. By Egorov's theorem, for any $\varepsilon >0$ and $k$ sufficiently large,
\begin{align*}
&\mathbb P( \mathcal W^{\bf g}(i)=0, \liminf_n\vert Z^{\bf g} _n(i, \cdot)\vert \ \geq 1)\\
&\leq 
\varepsilon + \mathbb P\left( \prod_{\ell=1}^p\mathbb P(\liminf_n\vert Z^{{\bf g}_k}_n(\ell, \cdot)\vert \ = 0)^{Z^{\bf g} _k(i,  \ell)}\geq 1-\varepsilon ;   \vert Z^{\bf g} _k(i, \cdot)\vert \ \geq 1\right)\\
 &\leq \varepsilon +{ 1\over   1-\varepsilon}
 \mathbb E\left[\prod_{\ell=1}^p\mathbb P(\liminf_n\vert Z^{{\bf g}_k}_n(\ell, \cdot)\vert \ = 0)^{Z ^{\bf g}_k(i,  \ell)};  \vert Z^{\bf g} _k(i, \cdot)\vert \ \geq 1\right] \\
 &=\varepsilon +{ 1\over   1-\varepsilon}
 \mathbb E\left[\mathbb P\left(\liminf_n\vert Z ^{\bf g}_n(i, \cdot)\vert \ = 0\Big\slash Z^{\bf g}_0, \ldots, Z^{\bf g}_n\right);  \vert Z^{\bf g} _k(i, \cdot)\vert \ \geq 1\right] \\
 &=\varepsilon +{ 1\over   1-\varepsilon}
 \mathbb P\left(\liminf_n\vert Z^{\bf g} _n(i, \cdot)\vert \ = 0;  \vert Z^{\bf g} _k(i, \cdot)\vert \ \geq 1\right). 
\end{align*}
Letting $k \to +\infty$, we obtain that $\displaystyle \mathbb P( \mathcal W^{\bf g}(i)=0, \liminf_n\vert Z^{\bf g} _n(i, \cdot)\vert \ \geq 1)\leq \varepsilon$ and the claim follows with $\varepsilon \to 0$.

\rightline{$\Box$}

  {\bf Proof of (\ref{convexg0})} We denote $\vert \cdot \vert _2$ the Euclidean norm on $\mathbb R^p$. The second inequality   is classical:    \[
  \vert  1-g_0^{(i)} (t)\vert\leq  \left\langle  \left(\nabla g_0^{(i)}\right)({\bf 1}), {\bf 1}-t\right\rangle = \langle  M_{g_0}(i, \cdot), {\bf 1}-t\rangle  \leq   \vert {\bf 1}-t\vert_2
  \preceq \vert M_{g_0}\vert \vert {\bf 1}-t\vert.  
  \]
 To prove the first inequality, let $[t', {\bf 1}]$ be the intersection of the cube $[0, 1[^p$ with the line passing through $t$ and ${\bf 1}$. By convexity of  $g_0^{(i)}$ on the segment $[t', {\bf 1}]$,  it holds that
  \[{1-g_0^{(i)}(t)\over \vert {\bf 1}-t\vert}\succeq {1-g_0^{(i)}(t)\over \vert {\bf 1}-t\vert_2}\geq {1-g_0^{(i)}(t')\over \vert {\bf 1}-t'\vert_2}\geq {1-g_0^{(i)}(t')\over \sqrt{p}}.
  \]
   The point $t'=(t'_1, \ldots, t'_p)$ belongs to the boundary of the cube $[0, 1[^p$ and at least one of its entries, say    $t'_j$, equals $0$; hence,  recalling that $\xi_0(i, \cdot)$ is a  $\mathbb N^p$-valued random variable with  generating function $g_0^{(i)}$, then
\[
 1-g_0^{(i)}(t')\geq1-g_0^{(i)}({\bf 1} - e_j)= 1-\mathbb P(\xi_0(i, j)=0)= \mathbb P(\xi_0(i, j)\geq 1)\geq \varepsilon.
 \] This achieves the proof.

\rightline{$\Box$}

\section{On the random environment}\label{Ontherandomenvironment}

In this section, we present the  random environment that we use  and introduce  some  considerable  classical change of measure   and its main properties.

Why this change of measure? The  bright idea introduced to study critical  branching processes in random environment is to assume first  that  the random walk $S_n$ is greater than some  constant $-a$, then  let $a\to +\infty$  (see for instance  \cite{DGV} and references therein). On the intermediate probability space,  for almost all environment with respect to the new probability measure, the Galton-Watson processe we  consider is  in varying environment and becomes a super critical process; we may thus apply Proposition \ref{varyingenvi}     to each one of these environment (quenched version).

  Recall that ${\bf f} = (f_n)_{n \geq 0}$ is a sequence of i.i.d. random variables  with values in $\mathcal G$.
   \subsection{Construction of a new  probability measure $\widehat {\mathbb{P}}_{\tilde x, a}$ }

The $\widetilde{P}_+$-harmonic function $V$ on $\mathbb{X} \times \mathbb{R}^+$  gives rise to a Markov kernel $\widetilde P _+^V$ on $\mathbb{X} \times \mathbb{R}^+ $ defined  formally by: 
 \[\widetilde P _+^V \phi =\frac{1}{V}\widetilde{P}_+(V \phi )\]
for any bounded measurable function $\phi$ on $\mathbb{X} \times \mathbb{R}^+$.  By Proposition \ref{prop1pham}, there exists $A>0$ such that  $V(\tilde x, a)>0$ whenever $a > A$; thus, for any 
  $\tilde x \in \mathbb{X}$, $a>A$ and $n \ge 1$, 
\begin{eqnarray}
 (\widetilde{P}_+^V)^n \phi  (\tilde x, a) 
 &=&{1\over V   (\tilde x, a)}  \mathbb{E}_{\tilde x, a} \left[ (V\phi)  (X_n, S_n); m_n >0\right]\notag.
\end{eqnarray}

 We introduce a 
change of probability measure on the canonical path space $ ((\mathbb{X} \times \mathbb{R})^{\otimes \mathbb{N}}, \sigma (X_n, S_n: n \geq 0), \theta)$ $^($\footnote{$\theta$ denotes   the shift operator on $(\mathbb{X} \times \mathbb{R})^{\otimes \mathbb{N}}$  defined by $\theta \Bigl((x_k, s_k)_{k\geq 0}\Bigr)=  (x_{k+1}, s_{k+1})_{k\geq 0} $ for any $(x_k, s_k)_{k\geq 0}$  in $(\mathbb{X} \times \mathbb{R})^{\otimes \mathbb{N}}$}$^)$   of  the Markov chain $(X_n, S_n)_{n \geq 0}$  from $\mathbb{P}$ to  the measure $\widehat {\mathbb{P}}_{\tilde x, a}$   characterized by the property that 
\begin{equation}  \label{eqn3.4}
\widehat {\mathbb{E}}_{\tilde x, a}[\phi (X_0, S_0,\ldots,X_k, S_k)] = \frac{1}{V{(\tilde x, a)}}{\mathbb{E} _{\tilde x, a}}[\phi ({X_0,\ldots, S_k})V(X_k, {S_k});\,\,{m_k} > 0]  
\end{equation}
for any positive Borel function $\phi$ on $(\mathbb{X}\times \mathbb{R})^{k+1}$. By Proposition \ref{prop1pham} and Theorem \ref{theo1pham}  
 
\begin{eqnarray} \label{eqn3.7}
 \lim_{n \to {+\infty} } \mathbb{E}_{\tilde x, a}[ \phi (X_0, \ldots,  S_k)\vert  m_n >  0]  &=& {1\over V(\tilde x,  a)} \mathbb{E}_{\tilde x, a} [V(X_k,   S_k)\phi (X_0,   \ldots,     S_k) ;  m_k  > 0]\notag \\
&=&  { \widehat{\mathbb E}}_{\tilde x,  a}[\phi (X_0,   \ldots,   S_k)],  
\end{eqnarray}
which clarifies the interpretation of $\widehat {\mathbb{P}}_{\tilde x,  a}$ (see  \cite{LPP} section 3.2 for the details).

\noindent This probability may be extended to the whole $\sigma$-algebra
$\sigma(f_n, Z_n: n \geq 0)$ as follows; the extension is done in three steps:    

\noindent {{\bf  Step 1.}} the marginal distribution of $\widehat{\mathbb{P}}_{\tilde x, a}$ on $\sigma({X_n, S_n: n\ge 0})$ is  $\widehat{\mathbb{P}}_{\tilde x, a }$ characterized by the property  (\ref{eqn3.4});

\noindent {{\bf  Step 2.}} for any $n \geq 0$, the conditional distribution of $ (f_0, \cdots, f_n)$ under $\widehat{\mathbb{P}}_{\tilde x, a}$ given $X_0 = \tilde x_0 =\tilde x,   \ldots, X_n = \tilde x _n, S_0 =s_0=a,  \ldots, S_n = s_n$ equals  the one of  $ (f_0, \cdots, f_n)$ under $\widehat{\mathbb{P}}$; namely,   for any measurable sets $G_0, \ldots, G_n$ in $\mathcal G $
and  all $(\tilde x _i)_{0\leq i \leq n}$ and $ (s_i)_{{0\leq i \leq n}}$ 
\begin{eqnarray*}  
\widehat{\mathbb{P}}_{\tilde x, a} (f_k \in G_k, 0\le k \le n\Big\slash X_i=\tilde x_i, S_i=s_i, 0\leq i \leq n) \qquad \qquad \qquad \qquad \qquad  \notag \\
 \qquad \qquad \qquad = \mathbb{P} (f_k \in G_k, 0\le k \le n\Big\slash X_i =\tilde x_i, S_i({\tilde x}, 0)=s_i,  0\leq i \leq n).
\end{eqnarray*}

\noindent {{\bf  Step 3.}} the conditional distribution of $(Z_n)_{n\ge 0}$ under $\widehat{\mathbb{P}}_{\tilde x, a}$ given ${\bf f}=(f_0 , f_1,\ldots ) $ is the same as under $\mathbb{P}$,;  namely, for any $n \geq 0$ and $1\leq i\leq p$, 
\begin{align*} 
\widehat{\mathbb{E}}_{\tilde x, a} \left[ s^{Z_n(i, \cdot)} \Big\slash Z_0, \ldots, Z_{n-1},f_0^{(i)}, f_1, \ldots, f_{n-1}\right]
&= \mathbb{E} \left[s^{Z_n (i, \cdot)}\Big\slash Z_0, \ldots, Z_{n-1}, f_0^{(i)}, f_1, \ldots, f_{n-1}\right]\notag
\\
&=  f_{n-1}(s)^{Z_{n-1}(i, \cdot)}.
\end{align*}

\subsection{ Some properties of the  probability measures  $\widehat{\mathbb P}_{\tilde x, a}, \tilde x \in \mathbb X , a \geq 0$}

 The following lemma extends property (\ref{eqn3.7}) to the $\sigma$-algebra $\mathcal F_\infty= \sigma (\vee_{k \geq 0} \mathcal F_k)$ where  $\mathcal F_k := \sigma \{f_\ell,  Z_\ell \ \vert \ 0\leq \ell \leq k \}$  for any $k\geq 0$.
 	
\begin{lem} \label{lemconvergence} 
	Assume that hypotheses { \rm H1--H5} hold for some $\delta >0$. Let $(Y_k)_{k \geq 0}$  be a sequence of  bounded real-valued random variables adapted to the filtration  $(\mathcal F_k)_{k \geq 0}$.
	\begin{enumerate}
		\item \cite{DV} For any $\tilde x \in  \X$ and $a>A$,
		\begin{equation}\label{probahatFk}
		\lim_{n \to +\infty } \E_{\tilde x, a} \big[ Y_k \ \vert \ \tau > n\big] = \widehat \E_{\tilde x, a}[Y_k].
		\end{equation}
		
		\item Moreover, if $(Y_k)_{k \geq 0}$ converges  in $\mathbb L^1(\widehat \P _{\tilde x, a})$ to some  random variable   $Y_\infty$,
		\begin{equation*} 
		\lim_{n \to +\infty } \E_{\tilde x, a} \big[Y_n \ \vert \ \tau>n\big] = \widehat \E_{\tilde x, a}[Y_\infty].
		\end{equation*}
		
	\end{enumerate}

\end{lem}

 \noindent {\bf Proof.}  Property (\ref{probahatFk}) is proved in \cite{DV}.  The second assertion  has an analogue version  in \cite{DV}, where  the almost-sure convergence is required;   in fact, the convergence in $L^1$ and the boundedness of the $Y_k$  suffice.
 
	For any $k \in \N$,  
	\begin{align*}
	 	\sqrt{n} \E_{\tilde x, a}[Y_n, \tau > n] 
	&=  \sqrt{n} \E_{\tilde x, a}[Y_k,  \tau > n] + \sqrt{n} \E_{\tilde x, a}[Y_n - Y_k,  \tau > n],  \notag 
	\end{align*}
	with
	\begin{eqnarray*}
	\lim_{n \to +\infty}\sqrt{n} \E_{\tilde x, a} [Y_k,  \tau > n] &=& \lim_{n \to +\infty}\sqrt{n} \E_{\tilde x, a} \big[Y_k \ \vert \ \tau > n\big] \P_{\tilde x, a}(\tau > n) \\
	&=&\frac{2}{\sigma \sqrt{2 \pi}} V(\tilde x, a) \widehat{\E}_{\tilde x, a} [Y_k],
	\end{eqnarray*}
by   (\ref{probahatFk}) and Theorem \ref{theo1pham}.
Since  $( \widehat{\E}_{\tilde x, a}[Y_k])_{k \geq 0}$ converges to $\widehat{\E}_{\tilde x, a}[Y_\infty]$ as $k \to +\infty$, it remains to prove that
\begin{eqnarray} \label{eqnYnk}
\lim_{k \to +\infty} \lim_{n \to {+\infty}} \sqrt{n} \E_{\tilde x, a}\big[\vert Y_n - Y_k\vert;  \tau > n\big] = 0.
\end{eqnarray}
 
We fix  $\rho > 1$ and decompose $\E_{\tilde x, a}\big[\vert  Y_n - Y_k\vert  ,  \tau > n\big]$ as
\begin{equation}\label{arkejh}
\E_{\tilde x, a}\big[\vert  Y_n - Y_k\vert  ,  \tau > n\big]= 
 \E_{\tilde x, a}\big[\vert  Y_n - Y_k\vert  , n < \tau < \rho n\big] + \E_{\tilde x, a}\big[\vert  Y_n - Y_k\vert  ,  \tau > \rho n\big].
\end{equation}
 For the first term in (\ref{arkejh}), since the random variables $Y_n$ are bounded, it is clear that 
 \begin{align*}
 \E_{\tilde x, a}\big[\vert  Y_n - Y_k\vert  , n < \tau < \rho n\big] &\preceq   \P _{\tilde x, a} (n < \tau < \rho n) \\
 & =   \P _{\tilde x, a} ( \tau > n) - \P _{\tilde x, a} ( \tau > \rho n). 
 \end{align*}
Therefore,  by Theorem \ref{theo1pham},  for any $k$ and $\rho >1$,
\begin{align*}
  \limsup_{n \to + \infty} \sqrt n  \E_{\tilde x, a} \big[ \vert  Y_n - Y_k\vert  , n < \tau < \rho n\big]   
&\preceq  \lim_{n \to + \infty} \sqrt n \P _{\tilde x, a} ( \tau > n) - \lim_{n \to + \infty} \sqrt n \P _{\tilde x, a} ( \tau > \rho n)    \\
&=  {2\over \sigma \sqrt{2\pi}  }  V(\tilde x, a)  \left( 1 - \frac{1}{\sqrt \rho}\right) \longrightarrow  0 \ \text{as} \ \rho \to 1.
\end{align*}
For the second term in (\ref{arkejh}),  we write
\begin{align*}
\E_{\tilde x, a}\big[ \vert  Y_n - Y_k\vert  ,  \tau > \rho n \big] 
&= \E_{\tilde x, a} \left[ \E \big[ \vert  Y_n - Y_k\vert  ,  \tau > \rho n \ \Big\slash \ \mathcal F_n\big] \right]  \\
&= \E_{\tilde x, a} \big[ \vert  Y_n - Y_k\vert   {\bf m}_{\rho n -n} (X_n, S_n),  \tau > n \big]  \\
&\preceq { 1\over   \sqrt{ n(\rho-1)}} \E_{\tilde x, a} \big[\vert  Y_n - Y_k\vert  V(X_n, S_n);   \tau > n \big] \\
&=   {1 \over \sqrt{n (\rho - 1)}} V(\tilde x, a)\widehat \E_{\tilde x, a} \big[ \vert  Y_n - Y_k\vert  \big]. 
\end{align*}
Hence, since $Y_n \to Y_\infty$ in $\mathbb L^1 (\widehat \P_{\tilde x, a})$,  
\begin{align*}
 & \limsup_{k \to + \infty} \limsup_{n \to + \infty} \sqrt n \E_{\tilde x, a} \big[ \vert  Y_n - Y_k\vert,  \tau > \rho n \big] \\
&\qquad \qquad \qquad \le  \frac{c }{\sqrt{\rho - 1}}  V(\tilde x, a)  \limsup_{k \to + \infty} \limsup_{n \to + \infty} \widehat \E_{\tilde x, a} \big[ \vert Y_n - Y_k\vert  \big]   \\
&\qquad \qquad \qquad =  \frac{c V(\tilde x, a)}{\sqrt{\rho - 1}}   \limsup_{k \to + \infty}  \widehat \E_{\tilde x, a} \big[ \vert  Y_\infty - Y_k\vert  \big]  =  0.
\end{align*}

\rightline{\rightline{$\Box$}}

The following statement plays a crucial role in the sequel. It was first proved in the multi-type context in \cite{LPP} (Lemma 3.1), when the generating functions are linear-fractional; then the general  case was  considered in  \cite{DV} (Lemma  7). We generalize these statements under weaker moment conditions.
\begin{lem} \label{serieswidetilde}Assume hypotheses { \rm H1--H6}   hold  for some $\delta>0$. Then, for any $\tilde x \in \mathbb X $ and $a>A$, 
\[
\sum_{n=0}^{+\infty}\widehat{\mathbb E}_{\tilde x, a} \left[ e^{-S_n}\right] <+\infty
\quad \text{and}\quad \sum_{n=0}^{+\infty}\widehat{\mathbb E}_{\tilde x, a} \left[ \eta_n e^{-S_n}\right] <+\infty. 
\]
\end{lem}
{\bf Proof.}   In order to ease the arguments for proving the first part of the statement, we begin by studying the second one. 
We fix $\tilde x \in \mathbb X , a>A$ and $n \geq 0$ and use Corollary  \ref{theo2pham2018} to control each term  $\widehat{\mathbb E}_{\tilde x, a}\left[\eta_{n } e^{-S_n}\right].$ By the definition of the probability measure  $\widehat{\mathbb P} _{\tilde x, a}$,
\begin{align*}
 \widehat{\mathbb E}_{\tilde x, a}\left[\eta_{n } e^{-S_n}\right] &=\widehat{\mathbb E}_{\tilde x, a}\left[ {\mu_{n }\over \vert M_{n }\vert^2}e^{-S_n}\right]\\
&\preceq
 \widehat{\mathbb E}_{\tilde x, a}\left[ {\mu_{n }\over \vert M_{n } X_{n }\vert^2}e^{-S_n}\right] \quad (\text{by Lemma} \  \ref{keylem})  \\
&\preceq
\int \widehat{\mathbb E}_{\tilde x, a}\left[ \mu_{n }\Big\slash X_n=\tilde y, S_n= s, X_{n+1 }=\tilde z, S_{n+1}= t\right] e^{s-2t} \\
&\qquad \qquad \qquad \qquad \widehat{\mathbb P}_{\tilde x, a} (X_n \in {\rm d} \tilde y, \ S_n\in{\rm d} s, X_{n+1} \in {\rm d} \tilde z, \ S_{n+1}\in{\rm d} t)\\
\\
&=
\int \mathbb E  \left[ \mu_{n}\Big\slash X_n=\tilde y, S_n= s, X_{n+1}=\tilde z, S_{n+1}= t\right] e^{s-2t} \\
&\qquad \qquad \qquad \qquad   \widehat{\mathbb P}_{\tilde x, a} (X_n \in {\rm d} \tilde y, \ S_n\in{\rm d} s, X_{n+1} \in {\rm d} \tilde z, \ S_{n+1}\in{\rm d} t).\\
\end{align*}
 Hence, by Proposition  \ref{prop1pham}
 \begin{align}\label{decompos}
 \widehat{\mathbb E}_{\tilde x, a}\left(\eta_{n } e^{-S_n}\right) 
&\preceq
\mathbb E_{\tilde x, a}\Bigl( \mathbb E  \left(\mu_{n }\Big\slash X_n, S_n, X_{n+1}, S_{n+1}\right)e^{S_n-2S_{n+1}} V(X_{n+1}, S_{n+1}); m_{n+1}>0
\Bigr)\notag
\\
&\preceq
\mathbb E_{\tilde x, a}\Bigl( \mathbb E  \left(\mu_{n }\Big\slash X_n, S_n, X_{n+1}, S_{n+1}\right)e^{S_n-2S_{n+1}} \vert S_{n+1}\vert  ; m_{n+1}>0
\Bigr)\notag
\\
&\preceq
\mathbb E_{\tilde x, a}\Bigl( \mathbb E  \left(\mu_{n }\Big\slash X_n, S_n, X_{n+1}, S_{n+1}\right)e^{S_n-2S_{n+1}}  \Bigl( \vert S_{n}\vert +\ln^+\vert M_{n }\vert\Bigr); m_{n+1}>0
\Bigr)
\end{align}
On the one hand
\begin{align}\label{omrvl}
\mathbb E_{\tilde x, a}&\Bigl( \mathbb E  \left(\mu_{n}\Big\slash X_n, S_n, X_{n+1}, S_{n+1}\right)e^{S_n-2S_{n+1}}  \vert S_{n}\vert  ; m_{n+1}>0
\Bigr)\notag
\\
&\preceq 
\mathbb E_{\tilde x, a}\Bigl( \mathbb E  \left(\mu_{n}\Big\slash X_n, S_n, X_{n+1}, S_{n+1}\right){e^{-S_n } \vert S_{n}\vert\over \vert M_{n}\vert ^2}  ; m_{n}>0
\Bigr)\notag\\
&= 
\mathbb E\Bigl(  {\mu_{n} \over  \vert M_{n}\vert ^2} e^{-S_n }   \vert S_{n}\vert ; m_{n}>0
\Bigr)\notag\\
&\leq \mathbb E_{\tilde x, a}\left(  {\mu_{n} \over  \vert M_{n}\vert ^2}\right)\times\mathbb E(e^{-S_n }  \vert S_{n} \vert ; m_{n}>0).
\end{align}
On the other hand  
\begin{align}\label{ozleinsf}
&
\mathbb E_{\tilde x, a}\Bigl( \mathbb E  \left(\mu_{n }\Big\slash X_n, S_n, X_{n+1}, S_{n+1}\right)e^{S_n-2S_{n+1}} \ln^+ \vert M_{n }\vert  ; m_{n+1}>0
\Bigr)\notag
\\
&\preceq 
\mathbb E_{\tilde x, a}\Bigl( \mathbb E  \left(\mu_{n }\Big\slash X_n, S_n, X_{n+1}, S_{n+1}\right)e^{-S_n }{ \ln^+ \vert M_{n }\vert\over \vert M_{n }\vert ^2}  ; m_{n}>0
\Bigr)
\notag
\\
&\leq \mathbb E\left(  {\mu_{n} \over  \vert M_{n}\vert ^2 } \ln^+ \vert M_{n}\vert\right)\times\mathbb E_{\tilde x, a}(e^{-S_n }; m_{n}>0).
\end{align}
By hypothesis H6, quantities $\displaystyle \mathbb E\left(  {\mu_{n} \over  \vert M_{n}\vert ^2}\right) $ and  $\displaystyle  \mathbb E\left(  {\mu_{n} \over  \vert M_{n}\vert ^2 } \ln^+ \vert M_{n}\vert\right) $ are both finite; furthermore,  Corollary \ref{theolocal} yields 
\[
n^{3/2}\mathbb E_{\tilde x, a}(e^{-S_n } \vert S_{n}\vert ; m_{n}>0)\preceq (1+a) \sum_{b \geq 0} (1+b)^2 e^{-b} <+\infty. 
\]
Finally, combining (\ref{decompos}), (\ref{omrvl}) and (\ref{ozleinsf}), we obtain that  
\[\displaystyle \sup_{n \geq 1} n^{3/2} \widehat{\mathbb E}_{\tilde x, a}\left(\eta_{n} e^{-S_n}\right) <+\infty
\] and the lemma follows.

\rightline{$\Box$}

As a direct consequence,  $\widehat{\mathbb P}_{\tilde x, a}$-almost surely, the environment ${\bf f}$ do satisfy the conclusions of Proposition \ref{varyingenvi}.  

\begin{coro} \label{corotildeP} Assume  that hypotheses { \rm H1--H6} hold  for some $\delta>0$. 
 
Then, 
for $\widehat{\mathbb P}_{\tilde x, a}$-almost all environment ${\bf f}= (f_n)_{n \geq 0}$,

1.   there exists a non-negative  random column vector ${\mathcal W}^{\bf f}= ({\mathcal W}^{\bf f}(i))_{1\leq i \leq p} $ such that for every $i, j \in \{1, \ldots, p\}$, as $n \to +\infty$, 
	\begin{equation} \label{w(i)bfrandom}
 {\mathcal W}^{\bf f}_n (i, j):= {Z^{\bf f}_n(i, j) \over\vert M_{0, n}^{\bf f}e_j\vert }\quad \xrightarrow{\mathbb L^2(\P _{{\bf f}}) } \quad  {\mathcal W}^{\bf f}(i).  
 \end{equation} 
 
If it is further assumed that there exists $\varepsilon, K>0$ such that all  the  generating functions $f_n, n\geq 0$, belong to $\mathcal G_{\varepsilon, K}$, then, 

2. the process $Z^{\bf f}$  becomes extinct  with probability  $ q^{\bf f}(i)<1$ for some (hence  every) $i\in \{1, \ldots, p\} $;

3. for any  $1\leq i\leq p,  $ it holds
\begin{equation*} 
 ({\mathcal W}^{\bf f}(i)>0) = \bigl(\bigcap_{n \geq 0} \vert  Z^{\bf f}_n(i, \cdot)\vert =0)\bigr)\quad \mathbb  P_{\bf f}\text{- a.s.}
\end{equation*}

\end{coro}

\noindent {\bf Proof.} By Lemma \ref{serieswidetilde},  for any $\tilde x \in \mathbb X$ and $a>A$, 
\[
\widehat{\mathbb E}_{\tilde x, a}\left[ \mathbb E\left[ \sum_{n=0}^{+\infty}  e^{-S_n}\Big\slash f_0, \ldots, f_{n-1}\right] \right] <+\infty
 \  \text{and} \  \widehat{\mathbb E}_{\tilde x, a}\left[ \mathbb E\left[ \sum_{n=0}^{+\infty} \eta_{n } e^{-S_n}\Big\slash f_0, \ldots, f_{n-1}\right] \right] <+\infty,\]
which yields, for $\widehat{\mathbb P}_{\tilde x, a}$-almost all ${\bf f}$,  
\[
\mathbb E\left[ \sum_{n=0}^{+\infty}  e^{-S_n(\tilde x, a)}\Big\slash { \bf f}\right] <+\infty
\quad \text{and}\quad \ \mathbb E\left[ \sum_{n=0}^{+\infty} \eta_{n } e^{-S_n(\tilde x, a)}\Big\slash  { \bf f}\right] <+\infty.\]
Hence, by Lemmas \ref{serieswidetilde} and  \ref{keylem}, for  $\widehat{\mathbb P}_{\tilde x, a}$-almost all ${\bf f}$ and  any $1\leq i, j \leq p$, 
 on the one hand,   
\[
\sum_{n=0}^{+\infty}  {1\over \vert M^{\bf f}_{0, n}(i, j)\vert} <+\infty,
 \] and on the other hand,   
\[
 \sum_{n=0}^{+\infty} { 1 \over \vert M_{0, n }^{\bf f}\vert  }\  {\bigl\vert B^{(i)}_{f_{n}} \bigr\vert
 \over \left\vert M_{f_{ n}} \right \vert  } \leq   \sum_{n=0}^{+\infty} { 1 \over \vert M_{0, n }^{\bf f}\vert  }\  \eta_{f_{n}} <+\infty.
 \] 
 Hence, $\widehat{\mathbb P}_{\tilde x, a}$-almost all environment ${\bf f}$ satisfy the hypotheses of Proposition \ref{varyingenvi},  Corollary \ref{corotildeP} follows immediately.
 
 \rightline{$\Box$}
 
 \subsection{On the  extinction  of $(Z_n(\tilde z, \cdot))_{n \geq 0}$ in random environment}\label{Ontheextinctionof}

The following result extends property  (\ref{extinctionAOP})  to  Galton-Watson  processes $(Z_n(\tilde z, \cdot))_{n \geq 0}$ with any initial population  $\tilde   z \in \mathbb N^p\setminus \{\tilde {\bf 0}\} $.

 Recall that for any $n \geq 0$ and $i \in \{1, \ldots, p\}, $ the probability of extinction at time $n$ of  $(Z_n(i, \cdot))_{n \geq 0}$, given the environment $\bf f$ (or equivalently given  $f_0, \ldots, f_{n-1}$) equals
\[q_{n, i}^{\bf f}  = \P (|Z_n(\tilde e_i, \cdot)| > 0 \ \slash  \ f_0, \ldots, f_{n-1}) = 1 - f_0^{(i)} f_1 \ldots f_{n-1} (\tilde {\bf 0}). 
\]
For any environment $\bf f$, the sequence $(q_{n, i}^{\bf f})_{n \geq 0}$ converges  to some limit, denoted $q_i^{\bf f}$.
Furthermore, by Corollary 5 in \cite{DV}, 
\begin{equation}\label{extinctni}
{1\over q_{n, i}^{\bf f}}={1\over 1- f_0^{(i)} f_1 \ldots f_{n-1} (\tilde {\bf 0})} \preceq { 1\over \vert M_{0, n}\vert}+\sum_{k=0}^{n-1}{\eta_{k }\over \vert M_{0, k}\vert}.
\end{equation}
By the branching property, for any $\tilde z=(z_1, \ldots, z_p) \in \mathbb N^p\setminus \{\tilde {\bf 0}\}$, 
\begin{align}\label{extinctprobz}
q^{\bf f}_{n,\tilde z} &:=    \P (|Z_n(\tilde{z},   \cdot)| > 0 \ \slash f_0, \ldots, f_{n-1}) \notag\\
&=  1 - \prod_{i=1}^p [ f_0^{(i)} f_1 \ldots f_{n-1} (\tilde {\bf 0}) ]^{z_i} 
=  1 - \prod_{i=1}^p [ 1 - q_{n, i}^{\bf f} ]^{z_i}.  
\end{align}
Let us denote $q^{\bf f}_{ \tilde z}$  the limit of the sequence $(q^{\bf f}_{n,\tilde z})_{n \geq 0}$. 

For any $\tilde  x \in \X$ and $a>A$, it holds that 
\[\P(|Z_n(\tilde{z},   \cdot)| > 0)= \E[ q_{n, \tilde z}^{\bf f}] \quad \text{and} \quad 
 \widehat{\P}_{\tilde x, a}(|Z_n(\tilde{z},   \cdot)| > 0)= \widehat{\E}_{\tilde x, a}[q_{n, \tilde z}^{\bf f}].
 \]
By the dominated convergence theorem,   
\[
\lim_{n \to +\infty}\P(|Z_n(\tilde{z},   \cdot)| > 0))=  \P(\cap_{n \geq 0}(|Z_n(\tilde{z},   \cdot)| > 0))) = \E[q_{ \tilde z}^{\bf f}]
\]
(resp. $\displaystyle 
\lim_{n \to +\infty}\widehat{\P}_{\tilde x, a}(|Z_n(\tilde{z},   \cdot)| > 0)
= \widehat{\E}[q_{ \tilde z}^{\bf f}]$).

These two limits are related to each other in the following way.

\begin{property} \label{prop5} Assume that   hypotheses { \rm H1--H6} hold for some $\delta>0$. Then for any $\tilde x \in \X$ and $\tilde z\in \mathbb N^p\setminus\{\tilde \bf 0\}$, 
	\begin{equation} \label{convpourz}
	\lim_{n \to +\infty} \sqrt{n} \P (\vert  Z_n (\tilde{z},   \cdot)\vert > 0 ) = \frac{2 }{\sigma \sqrt{2 \pi }}\lim_{a \to +\infty} V(\tilde x, a) \widehat{\E}_{\tilde x, a}(q_{\tilde z}^{\bf f}) \quad   =: \beta_{\tilde z} >0. 
	\end{equation}
	\end{property} 
  {\bf Proof. } We follow the   proof  detailed in  \cite{LPP}  and  \cite{DV} when $\tilde z=\tilde e_i$, it works along the same lines for general $\tilde z$. We fix $\tilde x \in \mathbb X$ and $a\geq 0$ and  decompose $\P (\vert  Z_n (\tilde{z},   \cdot)\vert > 0 )$ as 
\[
\underbrace{\P\left( \vert Z_n(\tilde{z},   \cdot )\vert > 0, \tau_{\tilde x, a} \le n \right)}_{A_n(\tilde x, a)}+
\underbrace{\P\left( \vert Z_n(\tilde{z},   \cdot )\vert > 0, \tau_{\tilde x, a} > n \right)}_{B_n(\tilde x, a)}.
\]
On the one hand,  by inequality (\ref{controlA}),  
\[
\limsup_{n \to +\infty} \sqrt{n} A_n(\tilde x, a)\preceq  \vert  z\vert  \     (1+a)e^{-a}  \longrightarrow \quad 0 \quad {\rm as} \quad   a\to +\infty.
\]
On the other hand,   it holds  $\P\left( \vert Z_n(\tilde{z},   \cdot )\vert > 0\slash \tau_{\tilde x, a} > n \right)=\mathbb E(q_{n, \tilde z}^{\bf f}\slash\tau_{\tilde x, a} > n )$; since $(q_{n, \tilde z}^{\bf f})_{n \geq 0} $ converges to $q_{  \tilde z}^{\bf f}$  in $\mathbb L^1(\widehat \P _{\tilde x, a})$,  Lemma \ref{lemconvergence}    (ii) yields
\[
\lim_{n \to +\infty} \P\left( \vert Z_n(\tilde{z},   \cdot )\vert > 0\slash \tau_{\tilde x, a} > n \right)=\widehat \E _{\tilde x, a}(q_{  \tilde z}^{\bf f}).
\] 
Hence, by using Theorem \ref{theo1pham}, we obtain, for any $\tilde x \in \mathbb X$ and $a\geq A$,
\begin{align*}
\lim_{n \to +\infty} B_n(\tilde x, a) &=\lim_{n \to +\infty} \P\left( \vert Z_n(\tilde{z},   \cdot )\vert > 0\slash \tau_{\tilde x, a} > n \right) \mathbb P( \tau_{\tilde x, a} > n)
\\
&=\frac{2 }{\sigma \sqrt{2 \pi }} V(\tilde x, a)  \widehat{\E}_{\tilde x, a}(q_{\tilde z}^{\bf f})<+\infty.
\end{align*}
Finally, 
\begin{align*}
\frac{2 }{\sigma \sqrt{2 \pi }}&V(\tilde x, a)  \widehat{\E}_{\tilde x, a}(q_{\tilde z}^{\bf f})
\leq \\
&\liminf_{n \to +\infty} \sqrt{n} \P (\vert  Z_n  (\tilde{z},   \cdot)\vert > 0 ) 
 \leq \limsup_{n \to +\infty} \sqrt{n} \P (\vert  Z_n  (\tilde{z},   \cdot)\vert > 0 )
 \\
 & \qquad \qquad \qquad \qquad \leq   c\  \vert  z\vert  \     (1+a)e^{-a}+\frac{2 }{\sigma \sqrt{2 \pi }}V(\tilde x, a)  \widehat{\E}_{\tilde x, a}(q_{\tilde z}^{\bf f})<+\infty.
\end{align*}
 In particular $\displaystyle \lim_{a \to +\infty} V(\tilde x, a)
\widehat{\E}_{\tilde x, a}(q_{\tilde z}^{\bf f})
$ exists and is finite; indeed,  the map $a \mapsto  V(\tilde x, a)  \widehat{\E}_{\tilde x, a}(q_{\tilde z}^{\bf f})$ is increasing (since  $a\mapsto   B_n(\tilde x, a) $ is also increasing) and bounded. Convergence (\ref{convpourz}) follows immediately and the limit $\beta_{\tilde z}$ is finite. 

It remains to prove that $\beta_{\tilde z}>0.$   Let $i_0\in \{1, \ldots, p\}$ such that $z_{i_0}\geq 1$;  by formula (\ref{extinctprobz}), it holds $q_{n, \tilde z}^{\bf f}\geq q_{n, i_0}^{\bf f}$ for all environment ${\bf f}$ so that $\widehat{\E}_{\tilde x, a}(q_{\tilde z}^{\bf f})\geq \widehat{\E}_{\tilde x, a}(q_{i_0}^{\bf f})$.  To conclude, it is sufficient to check that $\displaystyle \lim_{a \to +\infty} V(\tilde x, a)
\widehat{\E}_{\tilde x, a}(q_{i_0}^{\bf f})>0$; this is done in \cite{DV} and \cite{LPP}, and based on the following properties: 

(i)  the map $a \mapsto  V(\tilde x, a)  \widehat{\E}_{\tilde x, a}(q_{i_0}^{\bf f})$ is increasing;

(ii) $V(\tilde x, a)>0$ for $a\geq A$;

(iii)   $\displaystyle  \widehat{\E}_{\tilde x, a}\left({1\over q_{i_0}^{\bf f}}\right) \preceq  \widehat{\E}_{\tilde x, a} (\vert M_{0, n}\vert^{-1})+\sum_{k=0}^{n-1}\widehat{\E}_{\tilde x, a}\left({\eta_{k }\over \vert M_{0, k}\vert}\right)<+\infty$, hence $\widehat{\E}_{\tilde x, a} (  q_{i_0}^{\bf f})>0$.

\noindent (the property (iii) follows  using   formula 
  (\ref{extinctni}) and Lemma \ref{serieswidetilde}).

\rightline{$\Box$} 
 
Similarly, we need to extend property (\ref{w(i)bfrandom}) to   Galton-Watson  processes $(Z_n(\tilde z, \cdot))_{n \geq 0}$ with any initial population  $\tilde   z \in \mathbb N^p\setminus \{\tilde {\bf 0}\}$. The following statement is a direct consequence of  a combination of Corollary  \ref{corotildeP} and the branching property.

\begin{property} \label{propertyzW} Assume  that hypotheses  { \rm H1--H6} hold for some $\delta>0$. Then   
for  all $\tilde x \in \X, a>A$ and $\widehat{\mathbb P}_{\tilde x, a}$-almost all environment ${\bf f}= (f_n)_{n \geq 0}$,  any $ \tilde z\in \mathbb N^p  \setminus \{\bf 0\} $, and any $j \in \{1, \ldots, p\}$,  
\[
{\mathcal W}^{\bf f}_n (\tilde z, j):= {Z^{\bf f}_n(\tilde z, j) \over \vert M_{0, n}^{\bf f} e_ j\vert}\quad \xrightarrow{\mathbb L^2(\P _{{\bf f}}) } \quad 
  {\mathcal W}^{\bf f}(\tilde z):=  \sum_{i=1}^p \sum_{k=1}^{z_i}\mathcal W_k^{\bf f}(i), 
\]
where the random variables ${\mathcal W}^{\bf f}_{ k}$ for $ k \geq 0,$ are independent copies of $\mathcal{\mathcal W}^{\bf f}$.  

In particular, for any $j \in \{1, \ldots, p\}$, 
\[
\lim_{n \to +\infty}{\mathbb  E^{\bf f}[Z^{\bf f}_n(\tilde z, j)] \over \vert M_{0, n}^{\bf f} e_ j\vert} =
\lim_{n \to +\infty}{\sum_{i=1}^p z_i\mathbb  E^{\bf f} [Z^{\bf f}_n(i, j)] \over\vert M_{0, n}^{\bf f} e_ j\vert} = \langle z, \E [{\mathcal W}^{\bf f}] \rangle.
\]
If it is further assumed that there exist  $\varepsilon \in ]0,1[$ and $K>0$  such that $f_n \in G_{\varepsilon, K}$ for any $n \geq 0$, then 
\begin{equation}\label{nonextinctionrandomf}
({\mathcal W}^{\bf f}(\tilde z)>0) = \bigcap_{n \geq 0}( \vert  Z^{\bf f}_n(\tilde z, \cdot)\vert >0) \quad \mathbb  P_{\bf f}\text{- a.s.}
\end{equation}
\end{property}


\section{Proof of Theorem \ref{theo2}}\label{proofoftheorem2}

 By a standard argument in probability theory,  since the random variable  $Z_n(\tilde z, j)\slash \vert M_{0, n}e_j\vert$, for $ n \geq 0, \tilde z \in \N^p\setminus\{\tilde{\bf 0} \}$ and $1\leq j\leq p$,  are non-negative, it suffices  to prove that the  sequence of  Laplace transform 
 \[
  \lambda \mapsto   \E \left[ \exp{\left(-\lambda \frac{ Z_n(\tilde{z},  j)}{|M_0 \ldots M_{n-1} e_j|}\right)}\ \Big\slash \ \vert Z_n (\tilde{z},   \cdot)\vert  > 0  \right] 
\]
converges on $[0, +\infty[ $ to some  function  which is    continuous at $0$.

  We fix $z \in  \N^p\setminus\{\tilde{\bf 0} \}$ and $1\leq j \leq p$.  For any $\lambda \geq 0$, 
\begin{align*}
&\E \left[ \exp{-\left(\lambda\frac{- Z_n(\tilde{z},  j)}{|M_0 \ldots M_{n-1} e_j|}\right)} \ \Big\slash \  \vert Z_n (\tilde{z},   \cdot)\vert > 0  \right] \\
&\qquad \qquad \qquad\qquad \qquad \qquad
= \frac{\sqrt{n} \E \left[ \exp{ \left(-\lambda\frac{  Z_n(\tilde{z},  j)}{|M_0 \ldots M_{n-1} e_j|}\right)},  \vert Z_n (\tilde{z},   \cdot)\vert > 0  \right]   }{\sqrt{n} \P (\vert Z_n (\tilde{z},   \cdot)\vert > 0 )  }.
\end{align*}
By Property \ref{prop5}, it suffices to   prove that the sequence $(\phi_{n,  \tilde z, j})_{n \geq 0 }$ defined by
\[
\forall \lambda \geq 0, \quad \phi_{n,  \tilde z, j}(\lambda):= \sqrt{n} \E \left[ \exp{\left(-\lambda \frac{ Z_n(\tilde{z},  j)}{|M_0 \ldots M_{n-1} e_j|}\right)},  \vert Z_n (\tilde{z},   \cdot)\vert  > 0  \right] 
\]
converges  to some  function $\phi_{  \tilde z, j} : \ \mathbb R^+\to [0, 1]$ such that  \[ \lim_{\lambda \to 0^+}\phi_{\tilde z, j}(\lambda) = \phi_{\tilde z, j}(0) =   \beta_{\tilde z}.\]

A candidate for this limit is
\begin{equation*} 
 \phi_{\tilde z, j}(\lambda) =\frac{2}{\sigma \sqrt{2 \pi}} \sum_{k=0}^{+\infty} \E_{\tilde x, a}\left[  V(X_k, 0)  {\bf 1}_{(T_k =k)}  
\Psi(\lambda,  X_k, 0, Z_k(\tilde z, \cdot), {\bf f}\circ \theta^k)\right],
\end{equation*}
where
\begin{equation}\label{Psi}
\Psi(\lambda, \tilde x', a', \tilde z', {\bf g} ):= \widehat \E_{\tilde x', a'} \left[\exp\Bigl(- \lambda {\mathcal W^{{\bf g} }( \tilde{z}' )\over \langle \alpha, \E(\mathcal W^{{\bf g}})\rangle}\Bigr){\bf 1}_{\cap_{n>1}(\vert Z_{n}^{{\bf g} } ( \tilde{z}', \cdot)\vert > 0)} \right]
\end{equation}
for any $\lambda \geq 0, \tilde x'\in \X, a'>0, z' \in \N^{p}\setminus\{0\}$ and ${\bf g} \in \mathcal G^{\N}$.
  For any $n \geq 1$, we set $T_n=\max\{k \ \slash \ 0\leq k\leq n\ {\rm such \  that }\ S_k=m_n\}$; the random variable  $T_n$ satisfy the following simple properties:

-   $T_n \le n$  for any $n \geq 1$; 

-  $T_n$ does not depend on the value of $S_0$;

- let  $m_{k,n} := \min \{ S_{k+1} - S_k, \ldots, S_n - S_k \}$, then for any $0\leq k\leq n$, 
\[
(T_n = k) = (T_k = k) \cap (m_{k,n} >0).
 \]
These random variable  yields to the following decomposition
\begin{align*}
 \phi_{{ n, \tilde z, j}}(\lambda)&= \sqrt{n} \E \left[ \exp{\left(-\lambda \frac{ Z_n(\tilde{z},  j)}{|M_0 \ldots M_{n-1} e_j|}\right)};  \vert Z_n (\tilde{z},   \cdot)\vert > 0  \right] \notag \\
&= \sqrt{n} \E_{\tilde x, a}\left[ \exp{\left(-\lambda \frac{Z_n(\tilde{z},  j)}{|M_0 \ldots M_{n-1} e_j|}\right)};  \vert Z_n (\tilde{z},   \cdot)\vert > 0 , T_n = n \right] \notag \\
& \qquad \qquad  + \sqrt{n} \sum_{k=0}^{n-1} \E_{\tilde x, a}\left[ \exp{\left(-\lambda \frac{Z_n(\tilde{z},  j)}{|M_0 \ldots M_{n-1} e_j|}\right)};  \vert Z_n (\tilde{z},   \cdot)\vert > 0 , T_n = k \right] \notag \\
&= \underbrace{\sqrt{n} \E_{\tilde x, a}\left[ \exp{\left(-\lambda \frac{Z_n(\tilde{z},  j)}{|M_0 \ldots M_{n-1} e_j|}\right)};  \vert Z_n (\tilde{z},   \cdot)\vert > 0 , T_n = n \right] }_{\Sigma_1(n, \lambda)} \notag \\
& \qquad \qquad  +\underbrace{ \sqrt{n} \sum_{k=0}^{n-1} \E_{\tilde x, a}\left[ \exp{\left(-\lambda \frac{Z_n(\tilde{z},  j)}{|M_0 \ldots M_{n-1} e_j|}\right)};  \vert Z_n (\tilde{z},   \cdot)\vert > 0 , T_k = k, m_{k, n} > 0 \right] }_{\Sigma_2(n, \lambda)}.\notag
\end{align*}
The following lemma shows that  
\begin{equation}\label{sigma1n}
\lim_{n \to +\infty}\Sigma_1(n, \lambda)= 0, \quad \text{uniformly in } \quad \lambda \geq 0.
\end{equation}

\begin{lem} \label{lem4.2}
	There exists a positive constant $c$ such that for any $n \ge 1,    \tilde x \in  \X, a > 0$ and $  z \in \N ^p\setminus\{\bf 0\}$,
	\begin{eqnarray*}
	\P_{\tilde x, a}(\vert Z_n (\tilde{z},   \cdot)\vert > 0 , T_n = n ) \le c\frac{ \vert  z\vert }{n^{3/2}}.
	\end{eqnarray*}
\end{lem}
The term $\Sigma_2(n, \lambda)$ may  be decomposed as follows: for $1\leq \ell\leq n-1$ fixed, 
\begin{align*}
\Sigma_2(n, \lambda) &=  \underbrace{\sqrt{n} \sum_{k=0}^{\ell} \E_{\tilde x, a}\left[ \exp{\left(-\lambda \frac{Z_n(\tilde{z},  j)}{|M_0 \ldots M_{n-1} e_j|}\right)};  \vert Z_n (\tilde{z},   \cdot)\vert > 0 , T_k = k, m_{k, n} > 0 \right] }_{\Sigma_{2, 1}(n, \ell, \lambda)}\\
&  \qquad  + \underbrace{\sqrt{n} \sum_{k=\ell+1}^{n-1} \E_{\tilde x, a}\left[ \exp{\left(-\lambda \frac{Z_n(\tilde{z},  j)}{|M_0 \ldots M_{n-1} e_j|}\right)};  \vert Z_n (\tilde{z},   \cdot)\vert > 0 , T_k = k, m_{k, n} > 0 \right] }_{\Sigma_{2, 2}(n, \ell, \lambda)} 
\end{align*}
and we study separately  the two terms $\Sigma_{2,1}(n,\ell)$ and $\Sigma_{2,2}(n,\ell)$. Firstly, 
\begin{align*}
&\Sigma_{2,1}(n, \ell, \lambda) \\
&\quad =  \sum_{k=0}^{\ell} \frac{\sqrt{n}}{\sqrt{n - k}} \E_{\tilde x, a}\left[ \sqrt{n - k} \exp{\left(-\lambda \frac{Z_n(\tilde{z},  j)}{|M_0 \ldots M_{n-1} e_j|}\right)};  \vert Z_n (\tilde{z},   \cdot)\vert > 0 , T_k = k,m_{k, n} > 0 \right] \notag \\
&\quad = \sum_{k=0}^{\ell} \frac{\sqrt{n}}{\sqrt{n - k}}  \int df_0 \ldots df_{k-1} \   \delta_{M_{0, k}}(dM)  \ \int \P \big( Z_k(\tilde{z},   \cdot) \in dZ \ \vert \ f_0, \ldots, f_{k-1} \big)  {\bf 1}_{(T_k =k)}  {\bf 1}_{(\vert Z\vert  > 0)} \\
& \qquad \qquad \times \sqrt{n - k}\  \E_{\tilde x \cdot M, 0} \left[ \exp{\left(-\lambda \frac{Z^{{\bf f}\circ \theta^k}_{n-k} (Z, j)}{\vert M M_{0,n-k} \circ \theta^k e_j\vert }\right)}; \vert Z_{n-k}^{{\bf f}\circ \theta^k} (Z, \cdot)\vert > 0, \tau\circ \theta^k > n-k \right].
\end{align*}
Let us fix $0\leq k\leq \ell$ and set $Y_{n, k}(\lambda, Z, j):= \displaystyle \exp{\left(-\lambda \frac{Z_{n-k}^{{\bf f}\circ \theta^k} (Z, j)}{|M M_{0,n-k} \circ \theta^k e_j|}\right)}{\bf 1}_{(\vert Z_{n-k}^{{\bf f}\circ \theta^k} (Z, \cdot)\vert > 0)}$ for $n >k$. By Property \ref{propertyzW},  for  $\widehat{\mathbb P}_{\tilde x\cdot M, 0}$-almost all environment  $\bf g$, 

$\bullet$  the  sequence 
$ \displaystyle \left({Z^{\bf g}_{n-k} (Z, j)\over | M^{\bf g}_{0,n-k} e_j|}\right)_{n >k}$ converges  in    $\mathbb L^2(\P_{\bf g})$  to $\displaystyle {\mathcal W}^{\bf g}(Z)=\sum_{i=1}^p \sum_{l=1}^{Z_i}\mathcal W_l^{\bf g}(i) $ where the $W_l^{\bf g}$, for $ l \geq 1$ are independent copies of ${\mathcal W}^{\bf g}$; 

$\bullet$  let $\alpha_i :=\vert M(\cdot , i)\vert$  for $1\leq i\leq p$ and $\alpha= (\alpha_i)_{1\leq i\leq p}$, it holds that
\begin{align*}
{ | MM^{\bf g}_{0,n-k} e_j|\over | M^{\bf g}_{0,n-k} e_j|}&={  \sum_{i=1}^p\alpha_i  M^{\bf g}_{0,n-k} (i, j) \over \sum_{i=1}^p   M^{\bf g}_{0,n-k} (i, j)}\\
&={\sum_{i=1}^p \alpha_i\mathbb  E^{\bf g} [Z _{n-k}(i, j)] \over \mathbb  E^{\bf g}\vert Z _{n-k}(\cdot, j)\vert} \quad  \xrightarrow \quad \langle \alpha, \E [{\mathcal W}^{\bf g}]\rangle \quad \text{as} \quad   n \to +\infty;
\end{align*}

$\bullet\quad  \displaystyle  \lim_{n \to +\infty} {\bf 1}_{(\vert Z_{n-k}^{\bf g} (Z, \cdot)\vert > 0)}= {\bf 1}_{\cap_{n>k}(\vert Z_{n-k}^{\bf g}  (Z, \cdot)\vert > 0))}$.

 Hence,  for any $ 0\leq k \leq \ell$,    the sequences $\left(Y_{n, k}(\lambda, Z, j)\right)_{n >k} $  converge   in $\mathbb L^1(\widehat{\P}_{\tilde x\cdot M, 0})$   to the random variable
\[ Y_{\infty, k}(\lambda, Z):= \exp\left(- \lambda {{\mathcal W}^{{\bf f}\circ \theta^k}(Z)\over \langle \alpha, \E({\mathcal W}^{{\bf f}\circ \theta^k})\rangle}\right){\bf 1}_{\bigcap_{n>k}(\vert Z_{n-k}^{{\bf f}\circ \theta^k} (Z, \cdot)\vert > 0))}.
\]
 Lemma \ref{lemconvergence}  yields that
 \[
 \lim_{n \to +\infty} \sqrt{n - k}\  \E_{\tilde x \cdot M, 0} \left[ Y_{n, k}( \lambda)\Big\slash  \tau\circ \theta^k > n-k \right]= \widehat \E_{\tilde x \cdot M, 0} [Y_{\infty, k}( \lambda, Z)]. \]
Consequently
\begin{align*}
&\lim_{n \to +\infty}
\Sigma_{2,1}(n, \ell, \lambda) \\
& = \sum_{k=0}^{\ell} \int df_0 \ldots df_{k-1}  \   \delta_{M_{0, k}}(dM)  \ \int_{(\vert Z\vert >0)} \P( Z_k(\tilde{z},   \cdot) \in dZ | f_0, \ldots, f_{k-1})   {\bf 1}_{(T_k =k)}   \notag \\
& \qquad \qquad \times \frac{2}{\sigma \sqrt{2 \pi}} V(\tilde x \cdot M, 0)\   \widehat \E_{\tilde x \cdot M, 0} \left[ \exp\left(- \lambda {{\mathcal W}^{{\bf f}\circ \theta^k}(Z)\over \langle \alpha, \E({\mathcal W}^{{\bf f}\circ \theta^k})\rangle}\right){\bf 1}_{\cap_{n>k}(\vert Z_{n-k}^{{\bf f}\circ \theta^k} (Z, \cdot)\vert > 0)}\right] \\
&=\frac{2}{\sigma \sqrt{2 \pi}} \sum_{k=0}^{\ell} \E_{\tilde x, a}\left[   V(X_k, 0)   {\bf 1}_{(T_k =k)} 
\Psi(\lambda, X_k, 0, Z_k(\tilde z, \cdot), {\bf f}\circ \theta^k)
  \right]
\end{align*}
where $\Psi$ is defined in (\ref{Psi}).
Notice that,  by Lemma \ref{lem4.2},  for any $k \geq 1$, 
\begin{align}\label{sigma21MAJORATION}
0&\leq  \E_{\tilde x, a}\left[  V(X_k, 0) 
\Psi(\lambda,  X_k, 0, Z_k(\tilde z, \cdot), {\bf f}\circ \theta^k)\right] \notag\\
&\qquad \qquad\qquad \qquad \qquad \preceq
  \P_{\tilde x, a} (T_k =k, \vert Z_k(\tilde{z},   \cdot) \vert > 0)    \preceq    { \vert   z\vert \over k^{3/2}} 
\end{align}
 so that uniformly  in $  \lambda \geq 0$, 
\begin{equation}\label{sigma21nell}
 \lim_{\ell \to +\infty}\lim_{n \to +\infty}
\Sigma_{2,1}(n,\ell, \lambda)=\frac{2}{\sigma \sqrt{2 \pi}} \sum_{k=0}^{+\infty} \E_{\tilde x, a}\left[ V(X_k, 0) 
\Psi(\lambda, X_k, 0, Z_k(\tilde z, \cdot), {\bf f}\circ \theta^k)\right]
\end{equation}
 exists and is finite.
 
Let us control the term  $\Sigma_{2,2}(n, \lambda)$.
\begin{align*} 
&\Sigma_{2,2}(n, \lambda) \notag \\
 &= \sqrt{n} \sum_{k=l+1}^{n-1}  \int df_0 \ldots df_{k-1} \   \delta_{M_{0, k}}(dM)  \int \P( Z_k(\tilde{z},   \cdot) \in dZ | f_0, \ldots, f_{k-1})   {\bf 1}_{(T_k =k)}  {\bf 1}_{(\vert Z\vert  > 0)} \notag \\
& \qquad \qquad\qquad  \times  \E_{\tilde x \cdot M, 0} \left[ \exp{\left(-\lambda \frac{Z_{n-k} (Z, j)}{|M M_{0,n-k} e_j|}\right)}; \vert Z_{n-k} (Z, \cdot)\vert > 0 , \tau > n-k \right] \notag \\
& \le \sqrt{n} \sum_{k=l+1}^{n-1}  \int df_0 \ldots df_{k-1}  \delta_{M_{0, k}}(dM)\\
&\qquad \qquad \qquad \int \P( Z_k(\tilde{z},   \cdot) \in dZ | f_0, \ldots, f_{k-1}) {\bf 1}_{(T_k =k)} {\bf 1}_{(\vert Z\vert  > 0)} \P_{\tilde x \cdot M, 0} (\tau \circ \theta^k> n-k ) \notag \\
&= \sqrt{n} \sum_{k=l+1}^{n-1}  \E_{\tilde x, a}\left[ \P_{X_k, 0} (\tau\circ \theta^k > n-k) ; T_k = k; \vert Z_k(\tilde{z},   \cdot)\vert >0 \right]. \notag 
\end{align*}
By   Theorem \ref{theo1pham}  and Proposition \ref{prop1pham},  
\begin{align*}
\Sigma_{2,2}(n, \lambda)  & \preceq  \sum_{k=l+1}^{n-1} \frac{\sqrt{n}}{\sqrt{n -k}}  \P_{\tilde x, a} (T_k = k; \vert Z_k(\tilde{z},   \cdot)\vert >0) \notag \\
& \le \vert  z\vert \sum_{k=l+1}^{n-1} \frac{   \sqrt n}{\sqrt{n -k} \ \ k^{3/2}} \\
 &  \preceq \vert   z\vert   \left( \frac{1}{\sqrt l} + \frac{1}{\sqrt n}\right)  
\end{align*}
which readily implies, uniformly in $\lambda \geq 0$,
\begin{equation}\label{sigma22n}
\displaystyle \limsup_{l \to +\infty} \limsup_{n \to +\infty} \Sigma_{2,2}(n, \lambda) = 0
\end{equation}
We conclude by combining (\ref{sigma1n}), (\ref{sigma21nell}) and (\ref{sigma22n}). In particular, since the above convergences are uniform in $\lambda \geq 0$, it holds that $\displaystyle \lim_{\lambda \to 0^+}\phi_{\tilde z, j}(\lambda) = \phi_{\tilde z, j}(0)=\beta_{\tilde z}$.

Finally, let us prove  that $\nu_{\tilde z,   j}(\{0\})=0$  when the offspring generating functions belong to $\mathcal G_{\varepsilon, K}$. It suffices to prove that the Laplace transform of $\nu_{\tilde z,   j}$   (or equivalently the function  $\phi_{\tilde z, j}$) tends to $0$ as $\lambda \to +\infty$. Indeed,    by  (\ref{nonextinctionrandomf}), $\widehat{\P}_{\tilde x, a}$-almost surely, 
\begin{eqnarray*}
\Psi(\lambda,  X_k, 0, Z_k(\tilde z, \cdot), {\bf f}\circ \theta^k)&=& \widehat \E_{\tilde X_k, 0} \left[ \exp\Bigl(- \lambda {{\mathcal W}^{{\bf g} }( \tilde{z} )\over \langle \alpha, \E({\mathcal W}^{{\bf f\circ \theta^k}})\rangle}\Bigr){\bf 1}_{ ({\mathcal W}^{\bf f}(\tilde z)>0)}\right] \Big\slash _{\tilde z= Z_k(\tilde x, \cdot)}\\
&\longrightarrow & 0 \quad \text{as} \quad \lambda \to +\infty.
\end{eqnarray*}
Hence, 
by combining the Lebesgue dominated convergence theorem and  (\ref{sigma21MAJORATION}), 
\[
\lim_{\lambda\to +\infty} \phi_{ \tilde z, j}(\lambda) =0.\]
This achieves the proof.

\rightline{$\Box$}

\noindent It remains to prove Lemma \ref{lem4.2}.

\noindent {\bf  Proof of Lemma \ref{lem4.2}. }
 By the branching property,
for any  $\tilde z\in \N^p\setminus\{\tilde{\bf 0} \}$ and  $\mathbb P$-almost all  environment  $\bf f$, 
\[Z^{\bf f}_n(\tilde z, \cdot) = \sum_{i=1}^p \sum_{k=z_1+\ldots +z_{i-1} +1 }^{z_1+\ldots +z_{i}}Z_{n, k}^{\bf f} (i, \cdot),
\]
where the $Z_{n, k}^{\bf f}$, for $ k \geq 1, $ are i.i.d. copies of $Z_n^{\bf f}$; in particular, if  $\vert Z^{\bf f}_n(\tilde z, \cdot)\vert  >0$, then there exist $i$ and $k$ such that  $1\leq i\leq p$ and $z_1+\ldots +z_{i-1} +1\leq k\leq z_1+\ldots +z_{i}$ and  $\vert Z_{n, k}^{\bf f} (i, \cdot)\vert >0$. Hence,   noticing  that $T_n$ does not depend on the value of $(X_0, S_0)$  and using Lemma \ref{keylem},  for any $\tilde x \in \X$ and $a>A$,  
\begin{align*}
\P_{\tilde x, a}(\vert Z_n (\tilde{z}, &  \cdot)\vert > 0 , T_n = n )=\P_{\tilde x, 0}(\vert Z_n (\tilde{z},   \cdot)\vert > 0 , T_n = n ) \\
  &  \leq \sum_{i=1}^p z_i \mathbb P_{\tilde x, 0}(\vert Z_n(i, \cdot)\vert >0, T_n=n) \\
  & \leq \sum_{i=1}^p z_i \mathbb E_{\tilde x, 0} \big[ \vert Z_n(i, \cdot)\vert; T_n=n \big] \\
   & = \sum_{i=1}^p z_i \mathbb E_{\tilde x, 0} \big[ \mathbb E[\vert Z_n(i, \cdot)\vert\slash f_0, \ldots, f_{n-1}]; T_n=n \big] \\
 &=  \sum_{i=1}^p z_i \E_{\tilde x, 0} \big[ \vert   M_{0, n} \vert; T_n = n \big]\\
 & = \vert z\vert \E_{\tilde x, 0} \big[ \vert   M_{0, n} \vert; T_n = n \big]  \\
 &= \ \vert z\vert\  \E_{\tilde x, 0}  \big[ \vert   M_{0, n} \vert; S_n \le S_0, S_n \le S_1, \ldots, S_n \le S_{n-1} \big] \\
 &\le  \ \vert z\vert\ \E \big[ \vert   M_{0, n} \vert;\vert M_{0, n}\vert  \le c, \vert M_{0, n}\vert \le c\vert   M_{0,1}\vert \ldots, \vert M_{0, n}\vert \le c\vert   M_{0,n-1}\vert \big] \\
 &= \ \vert z\vert\ \E \big[  \vert M_{n, 0}\vert; \vert M_{n, 0}\vert \le c, \vert M_{n, 0}\vert \le c\vert M_{n,n-1}\vert, \ldots, \vert M_{n, 0}\vert \le c\vert M_{n,1} \vert \big]
\\
 & \qquad \qquad \text{since} \ (M_0, \ldots, M_{n-1}){\stackrel{\rm dist.}{=}}(M_{n-1}, \ldots, M_{0})\\
  \end{align*}
\begin{align*}
  \qquad \qquad   \qquad&\le  \ \vert z\vert\ \E \big[  \vert M_{n, 0}\vert; \vert M_{n, 0}\vert \le c, \vert M_{n-1,0}\vert \le c^2, \ldots, \vert M_{1,0}\vert \le c^2 \big] \\
 &\le c\ \vert z\vert\ \E \big[  \vert M_{n, 0} x\vert; \vert M_{n, 0}x\vert \le c^2,\vert M_{n-1,0}x\vert\leq c^2, \ldots, \vert M_{1,0}x\vert \le c^2 \big] \\
 &= c^3\ \vert z\vert\  \E \left[ \frac{1}{c^2} \vert M_{n, 0} x\vert; \frac{1}{c^2} \vert M_{n, 0}x\vert \le 1, \frac{1}{c^2} \vert M_{n-1,0} x\vert \le 1, \ldots, \frac{1}{c^2} \vert M_{1,0} x\vert \le 1 \right] \\
 &= c^3\ \vert z\vert\   \E_{x, - \ln c^2} \big[ \exp{(S'_n)}; S'_n \le 0; S'_{n-1} \le 0, \ldots, S'_1 \le 0 \big]
 \\
 & = 
 c^3\ \vert z\vert\   \E_{x, - \ln c^2} \big[ \exp{(S'_n)}; \tau'>n \big]
\end{align*}
with $S'_n=S_n'(  x, a) = a+\ln \vert M_{0, n}x\vert$ for any $\tilde x \in \X$ and $a\in \R$ and
$\tau'=\tau'_{x, a}=   \min \{ n \ge 1: S'_n (x, a) > 0 \}$.

Similar statements as  Theorem \ref{theo1pham}, Proposition \ref{prop1pham} and Corollary \ref{theolocal} also exist for the sequence $(S'_n(x, a))_{n \geq 0}$ and the stopping time $\tau'$; in particular,   there exists a positive constant  $c'$ such that for any $\tilde x \in  \mathbb X, a, b \in \R$ and $n \geq 1$,
\begin{equation*} \label{theolocal'}
	0\leq \mathbb P _ {x, a}(S'_n \in ]b-1, b ], \tau'  > n) \leq c' {(1+\vert a\vert )(1+\vert b\vert)\over n^{3/2} }.
\end{equation*} 
 Therefore, 
\begin{align*}
 \P_{\tilde x, a}(\vert Z_n (\tilde{z},  & \cdot)\vert > 0, T_n = n )\\
& \le c^3\ \vert z\vert\   \E_{x, - \ln c^2} \left[ \exp{(S'_n)}; \tau'>n \right]  \\
 & = c^3\ \vert z\vert\  \sum_{b\leq 0} e^{b} \ \P_{x, - \ln c^2}\left( S'_n \in ]b-1, b], \tau' > n \right)  \\
 & \leq  c^3 \ (1+\vert \ln c^2\vert) \ \ c' \ \vert z\vert\  \left(\sum_{b\leq 0} e^{b} (1+\vert b\vert)\right)  {1\over n^{3/2}}.\\ 
\end{align*}
\rightline{$\Box$}
\section{Proof of Proposition \ref{proposition2.3}}\label{ProofofProposition2.3}

We fix  $t \in \R, \tilde x \in X, a>A$ and $\tilde z \in \N^p\setminus\{\tilde{\bf 0} \}$.  By  Property  \ref{prop5}, we have to prove that the sequence 
\[
 \left( 
\sqrt{n}\P \left(  {S_n (\tilde x, a)\over \sqrt n} \le t\ ,   \vert Z_n(\tilde{z},   \cdot )\vert > 0 \right)\right)_{n \geq 0}
\] converges as $n \to +\infty$ and identify its limit.

For any   $b\geq 0, \rho \in ]0, 1[$ and $m\in \{1, \ldots, [\rho n]\}$,   we may decompose the quantity
\[ \sqrt{n}\P \left(  {S_n (\tilde x, a)\over \sqrt n} \le t\ ,   \vert Z_n(\tilde{z},   \cdot )\vert > 0 \right)\]
 as
\begin{align*}
& \underbrace{\sqrt{n}\P \left( \frac{S_n (\tilde x, a)}{\sqrt n} \le t, \vert Z_n(\tilde{z},   \cdot )\vert > 0, \tau_{\tilde x, b} \leq  n \right)}_{A_n(b)}+    \sqrt{n} \P \left( \frac{S_n (\tilde x, a)}{\sqrt n} \le t, \vert Z_n(\tilde{z},   \cdot )\vert > 0, \tau_{\tilde x, b} > n \right) \notag \\
&\qquad =  A_n(b)+  \sqrt{n}  \P \left( \frac{S_n (\tilde x, a)}{\sqrt n} \le t, \vert Z_{[\rho n]}(\tilde{z},   \cdot )\vert > 0, \tau_{\tilde x, b} > n \right)   \\
 &  \qquad \qquad \qquad  -  \underbrace{\sqrt{n}  \P \left( \frac{S_n (\tilde x, a)}{\sqrt n} \le t, \vert Z_{[\rho n]}(\tilde{z},   \cdot )\vert > 0, \vert Z_n(\tilde{z},   \cdot )\vert = 0, \tau_{\tilde x, b} > n \right)}_{B_n(b,\rho)}  \\
&\qquad =  A_n(b)-B_n(b,\rho)+  \underbrace{\sqrt{n}   \P \left( \frac{S_n (\tilde x, a)}{\sqrt n} \le t,\vert Z_m(\tilde{z},   \cdot )\vert > 0,  \tau_{\tilde x, b} > n \right)}_{C_n(b,  \rho, m)}\\
 &  \qquad \qquad \qquad  \qquad \qquad   
-  \underbrace{\sqrt{n} \P \left( \frac{S_n (\tilde x, a)}{\sqrt n} \le t, \vert Z_m(\tilde{z},   \cdot )\vert > 0, \vert Z_{[\rho n]}(\tilde{z},   \cdot )\vert = 0, \tau_{\tilde x, b} > n \right)  }_{D_n(b,  \rho, m)} \notag \\
  &\qquad = A_n(b)-B_n(b,\rho) +C_n(b,  \rho, m)-D_n(b,  \rho, m).
\end{align*}

We control these terms one by one.

 \noindent  {\bf Step 1. } {\it The sequence $(A_n(b))_{n\geq 0}$ converges to $A(b)\geq 0$  and $\displaystyle \lim_{b\to +\infty} A(b)=0.$}
 
  This is a direct consequence of the following inequality: for any $n \geq 1$ and $ b \geq 0$,
 \begin{equation}\label{controlA} 
	\sqrt{n}\P\left( \vert Z_n(\tilde{z},   \cdot )\vert > 0, \tau_{\tilde x, b} \le n \right) \le  c\  \vert  z\vert  \     (1+b)e^{-b},	\end{equation}
for some positive constant $c$.
Indeed, by  (\ref{generating})  and (Lemma \ref{keylem},  for all $\tilde x\in \X$ and $1\leq k\leq n$,  it holds $\mathbb P$-a.s.  that
\[
\P( \vert Z_n(i,   \cdot )\vert>0 \ \vert \ f_0, \cdots , f_{n-1} )=1-f_{0, n}^{(i)}(\tilde {\bf 0})\leq \left(\sum_{l=1}^n{1\over \vert M_{0, l} \vert}\right)^{-1}
\leq c\vert \tilde x M_{0, k}\vert
\]
so that 
\[\P( \vert Z_n(i,   \cdot )\vert>0 \ \vert \ f_0, \cdots , f_{n-1} )\leq  c e^{m_n(\tilde x, 0)}.   
\]
This yields
\begin{align*} 
 \sqrt n & \P \left( \vert Z_n(\tilde{z},   \cdot )\vert > 0, \tau_{\tilde x, b} \le n \right) \\
&\qquad \leq  \sqrt n \sum_{i =1}^p z_i\E  \left[ \P  \left(  \vert Z_n(i,   \cdot )\vert > 0\Big\slash f_0, \cdots , f_{n-1} \right), \tau_{\tilde x, b} \le n \right] \\  
&\qquad \leq  c\ \sqrt n \ \vert z\vert \ \E \left[  e^{m_n(\tilde x,  0)};\tau_{\tilde x, b} \le n \right]\\	
&\qquad \leq  c\ \sqrt n \ \vert z\vert \ \E \left[  e^{m_n(\tilde x,  0)};m_n(\tilde x,  0)<-b \right]\\ 
 &\qquad = c \sqrt n \ \vert z\vert \  \sum_{k=0}^{+\infty}e^{-k-b}\P  \left(-k-1-b\leq   m_n(\tilde x, 0) < -k-b\right)    \\
&\qquad \leq c \sqrt n \ \vert z\vert \  \sum_{k=0}^{+\infty}e^{-k-b}\P  \left(  m_n(\tilde x, 0) \geq -k-1-b  \right)\\
&\qquad = c \sqrt n \ \vert z\vert \  \sum_{k=0}^{+\infty}e^{-k-b}\P_{\tilde x, k+1+b} (  \tau>n)  \\
 &\qquad \preceq \vert z \vert \ e^{-b}  \sum_{k = 0}^{+\infty} (b+k+2) e^{-k} \preceq \vert z \vert \ (1+b) e^{-b}
 \quad \text{by Proposition} \ \ref{prop1pham}\ \text{and Theorem} \  \ref{theo1pham}. 
\end{align*}

\noindent {\bf Step 2. } {\it  For any $ b \geq 0, \rho \in ]0, 1[$ and $0\leq m\leq [\rho n]$,  the sequence  $(D_n(b, \rho, m))_{n \geq 0}$ converges to $0$.}

 It suffices to prove that
 \begin{equation}\label{step2Bis}
	\lim_{n \to +\infty } \sqrt{n} \P\left(\vert Z_m(\tilde{z},   \cdot )\vert > 0, \vert Z_{[\rho n]}(\tilde{z},   \cdot )\vert = 0,  \tau_{\tilde x, b} > n \right) = 0.
	\end{equation}
	For $1 \le m \le  [\rho n]$,
\begin{align*}
 \P_{\tilde x, b}&\left(\vert   Z_m(\tilde z ,   \cdot )\vert > 0, \vert Z_{[\rho n]}(\tilde z,   \cdot )\vert = 0,  \tau > n \right)  \\
&=  \P_{\tilde x, b}\left(\vert Z_m(\tilde{z},   \cdot )\vert > 0, \tau > n \right) - \P_{\tilde x, b}\left( \vert Z_{[\rho n]}(\tilde{z},   \cdot )\vert > 0,  \tau > n \right)   \\
&=  \E_{\tilde x, b}\Bigl[  \P (\vert Z_m(\tilde{z},   \cdot )\vert > 0 \ \vert \ f_0, \ldots , f_{n-1}) 
-   \P (\vert Z_{[\rho n]}(\tilde{z},   \cdot )\vert > 0 \ \vert \ f_0, \ldots , f_{[\rho n]-1}); \ \tau > n  )\Bigr]   \\
&=  \E_{\tilde x, b}\left[ q_{m, \tilde z}^{\bf f}  - q_{[\rho n], \tilde z}^{\bf f} ; \ \tau > n  \right].   \\
\end{align*}

Hence, by Theorem \ref{theo1pham}, 
\begin{align*}
  \P_{\tilde x, b}  (\vert Z_m(\tilde{z},   \cdot )\vert > 0,& \vert Z_{[\rho n]}(\tilde{z},   \cdot )\vert = 0,  \tau > n )     \\
&\le  c  \E_{\tilde x, b} \left[ ( q_{m, \tilde z}^{\bf f}  - q_{[\rho n], \tilde z}^{\bf f}   ) \frac{V(X_{[\rho n]}, S_{[\rho n]})}{\sqrt{n - [\rho n]}} ; \tau > [\rho n]  \right]    \\
&\le   {c\over \sqrt{n(1 - \rho) }} \E \left[ ( q_{m, \tilde z}^{\bf f}  - q_{[\rho n], \tilde z}^{\bf f} ) V(X_{[\rho n]}, S_{[\rho n]}); \tau > [\rho n] \right]     \\
&=   {c\over \sqrt{n(1 - \rho) }} V(\tilde x, a) \widehat \E _{\tilde x, b} \left[  q_{m, \tilde z}^{\bf f}  - q_{[\rho n], \tilde z}^{\bf f}   \right].     
\end{align*}
Therefore,
\begin{align*}
  \lim_{m \to +\infty} \lim_{n \to +\infty } \sqrt{n} \P_{\tilde x, b}& (\vert Z_m(\tilde{z},    \cdot )\vert > 0, \vert Z_{[\rho n]}(\tilde{z},   \cdot )\vert  = 0,  \tau > n )     \\
&\le  \frac{c}{\sqrt{1 - \rho }} V(\tilde x, a) \lim_{m \to +\infty} \underbrace{\lim_{n \to +\infty } \widehat \E _{\tilde x, b} [ q_{m, \tilde z}^{\bf f}  - q^{\bf f} _{[\rho n], \tilde z} ]}_{ \widehat \E _{\tilde x, b} [ q_{m, \tilde z}^{\bf f}  - q^{\bf f}_{\tilde z}  ] }=0.\end{align*}
 where the last equality  is a direct consequence of the preamble of subsection \ref{Ontheextinctionof}.

 \noindent  {\bf Step 3. } {\it For any $ b \geq 0$ and $\rho \in ]0, 1[$,  the sequence  $(B_n(b,\rho))_{n \geq 0}$ converges to $0$.}

We write
\begin{align*}
 0\leq B_n(b,\rho)&= \P \left( \frac{S_n (\tilde x, a)}{\sqrt n} \le t, \vert Z_{[\rho n]}(\tilde{z},   \cdot )\vert > 0, \vert Z_n(\tilde{z},   \cdot )\vert = 0, \tau_{\tilde x, b} > n \right)   \\
&  \le  \P \left(  \vert Z_{[\rho n]}(\tilde{z},   \cdot )\vert > 0, \vert Z_n(\tilde{z},   \cdot )\vert = 0, \tau_{\tilde x, b} > n \right)   \\
& =   \P_{\tilde x, b}\left(  \vert Z_{[\rho n]}(\tilde{z},   \cdot )\vert > 0,\tau > n \right) - \P_{\tilde x, b}\left(  \vert Z_n(\tilde{z},   \cdot )\vert > 0, \tau > n \right).   
\end{align*}
By Lemma \ref{lemconvergence}  and Theorem \ref{theo1pham} ,  
\begin{equation*}
\lim_{n \to +\infty } \sqrt{n} \P_{\tilde x, b}\left( \vert Z_n(\tilde{z},   \cdot )\vert > 0, \tau > n \right)  = \frac{2 }{\sigma \sqrt{2 \pi}}V(\tilde x, b) \widehat \E_{\tilde x, b}\left[ q_{\tilde{z}}^{\bf f}   \right]
\end{equation*}
and  it suffices to check that the sequence $\left( \sqrt{n} \P_{\tilde x, b}\left(  \vert Z_{[\rho n]}(\tilde{z},   \cdot )\vert > 0,\tau > n \right)\right)_{n \geq 0} $ converges to the same limit.
Indeed, for $1 \le m \le [\rho n]$,
\begin{align*}
  \sqrt{n} \P_{\tilde x, b}&\left(  \vert Z_{[\rho n]}(\tilde{z},   \cdot )\vert > 0,\tau > n \right)   \\
& = \sqrt{n}  \P_{\tilde x, b}\left( \vert Z_m (\tilde{z},   \cdot )\vert > 0, \tau > n \right) -\sqrt{n}  \P_{\tilde x, b}\left( \vert Z_m (\tilde{z},   \cdot )\vert > 0, \vert Z_{[\rho n]}(\tilde{z},   \cdot )\vert = 0,\tau > n \right)   
\end{align*}
with 

$\bullet\quad \displaystyle\lim_{n \to +\infty} \sqrt{n} \P_{\tilde x, b}\left( \vert Z_m (\tilde{z},   \cdot )\vert > 0, \tau > n \right) = \frac{2 V(\tilde x, b)}{\sigma \sqrt{2 \pi}} \widehat \P_{\tilde x, b}(\vert Z_m(\tilde{z},   \cdot )\vert > 0 )
$, by Lemma \ref{lemconvergence} and   Theorem \ref{theo1pham};

$\bullet  \quad   \displaystyle \lim_{n \to +\infty} \sqrt{n} \P_{\tilde x, b}\left( \vert Z_m (\tilde{z},   \cdot )\vert > 0, \vert Z_{[\rho n]}(\tilde{z},   \cdot )\vert = 0,\tau > n \right) = 0$ by (\ref{step2Bis}) of Step 2.

\noindent Hence
\[
\lim_{m \to +\infty} \lim_{n \to +\infty} \sqrt{n} \P_{\tilde x, b}\left( \vert Z_{[\rho n]} (\tilde{z},   \cdot )\vert > 0, \tau > n \right) = \frac{2 V(\tilde x, b)}{\sigma \sqrt{2 \pi}} \widehat \P_{\tilde x, b} \Big(  q_{\tilde{z}}^{\bf f}  \Big)
\]
and the proof is complete.

 \noindent {\bf  Step 4.}
	{\it For  any $b \geq 0$ and $\rho \in ]0, 1[,$}
	\begin{equation*}\label{controlC}
	 \lim_{m \to +\infty} \lim_{n \to +\infty } C_n(b, \rho, m)
	 = \frac{2}{\sigma \sqrt{2 \pi}} V(\tilde x, b)\  \widehat \E _{x, b} [q^{\bf f}_{\tilde z}  ] \ \Phi^+ \left( \frac{t}{\sigma} \right).  
	\end{equation*}

Assume that $n \ge 2m$. 
On the one hand, when $t < 0$, the quantity  $t + \frac{b-a}{\sqrt n}$ becomes negative when $n$ is great enough, in which case   $\left( \frac{S_n}{\sqrt n} \le t + \frac{b-a}{\sqrt n} \right)\cap (\tau>n) = \emptyset$;  therefore, the above limit holds in this case.
On the other hand, when $t\geq 0$,
\begin{align*}
&  \P \left( \frac{S_n (\tilde x, a) }{\sqrt n} \le t,\vert Z_m(\tilde{z},   \cdot )\vert > 0,  \tau_{\tilde x, b} > n \right)   \\
&\quad =  \P_{\tilde x, b} \left( \frac{S_n}{\sqrt n} \le t + \frac{b-a}{\sqrt n},\vert Z_m(\tilde{z},   \cdot )\vert > 0,  \tau > n \right)   \\
&\quad =  \E_{\tilde x, b} \left[ \P_{\tilde x, b} \left( \frac{S_n}{\sqrt n} \le t + \frac{b-a}{\sqrt n}, \vert Z_m(\tilde{z},   \cdot )\vert > 0,  \tau > n
\ \Big\slash \ f_0, \ldots, f_{m-1}, Z_0, \ldots, Z_m \right) \right]   \\
&\quad =  \E_{\tilde x, b} \left[ \vert Z_m(\tilde{z},   \cdot )\vert > 0, \tau > m ,\right.   \\
& \qquad  \left.   \P_{\tilde x, b} \left( \frac{S_m + S_{n-m} \circ \theta^m (X_m, S_m)}{\sqrt{n-m}} \le    t_{n, m},  \tau_{(X_m, S_m)} > n - m \ \Big\slash \ f_0, \ldots, f_{m-1}, Z_0, \ldots, Z_m \right) \right],   
\end{align*}
where $   t_{n, m} = \left( t + \frac{b-a}{\sqrt n} \right) \frac{\sqrt{n}}{\sqrt{n-m}}$.
By Corollary \ref{theo2pham2018}, as $n \to +\infty$,  
\[
\sqrt{n} \P_{\tilde x, b}\left( \frac{S_n}{\sqrt n} \le s, \tau > n \right) \to \frac{2 V(\tilde x, a)}{\sigma \sqrt{2 \pi}} \Phi ^+ \left( \frac{s}{\sigma} \right) =  \frac{2 V(\tilde x, a)}{\sigma \sqrt{2 \pi}} \left(1 - \exp{\left(-\frac{s^2}{2 \sigma ^2} \right)} \right),
\]
then 
\begin{align*}
  \sqrt{n -m}\  \P \Bigl(  &\frac{S_m (\tilde x, a) + S_{n-m} \circ   \theta^m (X_m, S_m)}{\sqrt{n-m}}  \leq     t_{n, m},   \tau_{(X_m, S_m)} > n - m\Big\slash  f_0, \ldots, f_{m-1}, Z_0, \ldots, Z_m \Bigr)   \\
&=  \frac{2}{\sigma \sqrt{2 \pi}} V(X_m, S_m) \left(1 - \exp{\left(-\frac{   t_{n, m} ^2}{2 \sigma ^2} \right)}\right)  (1 + o(n-m)).  
\end{align*}
 Therefore, 
\begin{align*}
& \sqrt{n}  \P \left( \frac{S_n (\tilde x, a)}{\sqrt n} \le t,\vert Z_m(\tilde{z},   \cdot )\vert > 0,  \tau_{\tilde x, b} > n \right)   \\
&=  \frac{2}{\sigma \sqrt{2 \pi}} \frac{ \sqrt{n}}{\sqrt{n -m}}\E_{\tilde x, b} \left[  V(X_m, S_m)(1 + o(n-m) ) \left(1 - \exp{\left(-\frac{   t_{n, m} ^2}{2 \sigma ^2} \right)} \right) , \vert Z_m(\tilde{z},   \cdot )\vert > 0, \tau > m \right]   \\
&=   
\frac{2}{\sigma \sqrt{2 \pi}} 
 \frac{ \sqrt{n}}{\sqrt{n -m}}
 V(\tilde x, b)\ 
\widehat \E_{\tilde x, b} \left[ 
\left(  1 - \exp{ \left(-\frac{   t_{n, m} ^2}{2 \sigma ^2} \right)}
 \right)   (1 + o(n-m)); \vert Z_m(\tilde z,   \cdot )\vert > 0  \right] \\
 &\qquad \longrightarrow
  \frac{2}{\sigma \sqrt{2 \pi}}   \left(  1 - \exp{\left(-\frac{   t^2}{2 \sigma ^2} \right)} \right)  V(\tilde x, b)\ \widehat \P_{\tilde x, b}( \vert Z_m(\tilde{z},   \cdot )\vert > 0 )\quad \text{as} \quad n \to +\infty.    \\
\end{align*}
Finally 
\begin{align*}
   \lim_{m \to +\infty} \lim_{n \to +\infty } \sqrt{n} &\P \left( \frac{S_n (\tilde x, a)}{\sqrt n} \le t,\vert Z_m(\tilde{z},   \cdot )\vert > 0,  \tau_{\tilde x, b} > n \right)   \\
&=  \frac{2}{\sigma \sqrt{2 \pi}}   \left(  1 - \exp{\left(-\frac{   t^2}{2 \sigma ^2} \right)}  \right) V(\tilde x, b)\  \widehat \P_{\tilde x, b} \left( \bigcap_{m \ge 0} [\vert Z_m(\tilde{z},   \cdot )\vert > 0] \right)   \\
&=  \frac{2}{\sigma \sqrt{2 \pi}}   \Phi^+ \left( \frac{t}{\sigma} \right) V(\tilde x, b)\ \widehat \P_{\tilde x, b} \left(q_{\tilde{z}}^{\bf f}  \right).     
\end{align*}

	 \noindent  {\bf Step 5.} {\it Conclusion}

By the four previous steps  and Property \ref{prop5}, letting $ n \to+\infty$, then $m\to +\infty$ and at last  $b\to +\infty$, we obtain that
\[ \lim_{n \to +\infty}\sqrt{n}\P \left(  {S_n (\tilde x, a)\over \sqrt n} \le t\ ,   \vert Z_n(\tilde{z},   \cdot )\vert > 0 \right)
= \frac{2}{\sigma \sqrt{2 \pi}} \Phi^+ \left( \frac{t}{\sigma} \right).  \]

\rightline{$\Box$}
\section{Proof of Theorem \ref{MAINTHEO}} \label{proofofmaintheo}

 Let $\varepsilon >0$. Then
\begin{align*} 
  \P \Bigl( \Big\vert \frac{\ln  \vert Z_n (\tilde z, \cdot)\vert}{\sqrt n} &- \frac{S_n (\tilde x, a)}{\sqrt n} \Big\vert \ge \varepsilon \ \Big\slash \ \vert Z_n (\tilde z, \cdot)\vert >0 \Bigr)   \\
&= \P \left( | \ln  \vert Z_n (\tilde z, \cdot)\vert - S_n (\tilde x, a) | \ge  \varepsilon {\sqrt n} \ \Big\slash \ \vert Z_n (\tilde z, \cdot)\vert >0 \right)   \\
&= \P \left( \frac{\vert Z_n (\tilde z, \cdot)\vert}{e^{S_n (\tilde x, a)}} \ge e^{ \varepsilon \sqrt n } \ \Big\slash \  \vert Z_n (\tilde z, \cdot)\vert >0 \right) 
\\
&
 \qquad \qquad \qquad \qquad    
 + \P \left( \frac{\vert Z_n (\tilde z, \cdot)\vert}{e^{S_n (\tilde x, a)}}  \le e^{- \varepsilon \sqrt n} \ \Big\slash \ \vert Z_n (\tilde z, \cdot)\vert >0 \right)\\
 &\leq \sum_{j=1}^p  
\P \left( \frac{  Z_n (\tilde z, j) }{e^{S_n (\tilde x, a)}} \ge e^{ \varepsilon \sqrt n }/p \ \Big\slash \ \vert Z_n (\tilde z, \cdot)\vert >0 \right) 
\\
& \qquad \qquad \qquad \qquad   + \P \left( \frac{  Z_n (\tilde z, 1) }{e^{S_n (\tilde x, a)}}  \le e^{- \varepsilon \sqrt n } \ \Big\slash \ \vert Z_n (\tilde z, \cdot)\vert >0 \right) 
\\
& \leq \sum_{j=1}^p  
\P \left( \frac{  Z_n (\tilde z, j) }{\vert M_{0, n}e_j\vert} \ge c \   e^{a+ \varepsilon \sqrt n }/p \ \Big\slash \ \vert Z_n (\tilde z, \cdot)\vert >0 \right)
\\
& \qquad \qquad \qquad \qquad    + \P \left( \frac{  Z_n (\tilde z, 1) }{\vert M_{0, n}e_1\vert}  \le {1\over c} \ e^{a- \varepsilon \sqrt n } \ \Big\slash \ \vert Z_n (\tilde z, \cdot)\vert >0 \right) 
\\
& \qquad \qquad \qquad \qquad   \text{(where  $c$ is the constant which appears in Lemma \ref{keylem})}
\\
&= \sum_{j=1}^p    \left(1-  \P \left( \frac{\vert Z_n (\tilde z, j)\vert}{\vert M_{0, n}e_j\vert}< c \   e^{a+ \varepsilon \sqrt n }/p \ \Big\slash \ \vert Z_n (\tilde z, \cdot)\vert >0 \right)\right)  
 \\
& \qquad \qquad \qquad \qquad   + \P \left( \frac{ Z_n (\tilde z,  1) }{\vert M_{0, n}e_1\vert}  \le  {1\over c} \ e^{a- \varepsilon \sqrt n }  \ \Big\slash \ \vert Z_n (\tilde z, \cdot)\vert >0 \right). 
\end{align*}
 Fix $A > 1$  then   there exists a number   $n_A$ great enough such that  $c\ e^{a+ \varepsilon \sqrt n_A }/p > A$ and ${1\over c} \ e^{a- \varepsilon \sqrt n_A } <1/A$ ;  without loss of generality, we assume $\nu_{\tilde z,  j}(\{{1\over A}, A\})= 0$. Hence, for any $n \geq n_A$,
 \begin{align*} 
  \P \Bigl( &\Big\vert \frac{\ln  \vert Z_n (\tilde z, \cdot)\vert}{\sqrt n}  - \frac{S_n (\tilde x, a)}{\sqrt n} \Big\vert \ge \varepsilon \ \Big\slash \ \vert Z_n (\tilde z, \cdot)\vert >0 \Bigr)   \\
&\leq   
 \sum_{j=1}^p    \left(1-  \P \left( \frac{\vert Z_n (\tilde z, j)\vert}{\vert M_{0, n}e_j\vert}< A \ \Big\slash \ \vert Z_n (\tilde z, \cdot)\vert >0 \right)\right)  
  + \P \left( \frac{ Z_n (\tilde z,  1) }{\vert M_{0, n}e_1\vert}  \le  {1\over A}  \ \Big\slash \ \vert Z_n (\tilde z, \cdot)\vert >0\right);  
\end{align*}
hence, by Theorem \ref{theo2},  
\[
 \limsup_{n \to +\infty} \P \Bigl( \Big\vert \frac{\ln  \vert Z_n (\tilde z,  \cdot )\vert}{\sqrt n} - \frac{S_n (\tilde x, a)}{\sqrt n} \Big\vert \ge \varepsilon \ \Big\slash \ \vert Z_n (\tilde z, \cdot)\vert >0 \Bigr)\leq   \sum_{j=1}^p  \left(1- \nu_{\tilde z,   j}([0, A ])\right)   +\nu_{\tilde z,   1}([0, 1/A]).
\]
Since the  $\nu_{\tilde z,   j}$ are probability measures on $]0, +\infty[$,  it holds $\nu_{\tilde z,   j}([0, A ]) \to 1 $ and $\nu_{\tilde z,  j}([0, 1/A]) \to 0 $ for any $1 \le j \le p$,  as $A\to +\infty$. This yields 
\[
 \lim_{n \to +\infty} 
  \P \Bigl( \Big\vert \frac{\ln  \vert Z_n (\tilde z, \cdot )\vert}{\sqrt n} - \frac{S_n (\tilde x, a)}{\sqrt n} \Big\slash \ge \varepsilon \ \Big\vert \ \vert Z_n (\tilde z, \cdot)\vert >0 \Bigr)  =0.
  \]
  We complete the proof by combining  Proposition \ref{proposition2.3}   and Slutsky's lemma.  
  
  \rightline{$\Box$}


\begin{thebibliography}{9}

 \bibitem{A}\textsc{ Afanasyev V.I. } (2001)
 A functional limit theorem for a critical branching process in a random environment, Discrete Math. Appl.,   \textbf{11},   \textbf{6},   587--606. 

\bibitem{ABKV}
\textsc{ Afanasyev V. I., B\"oinghoff C., Kersting G., \& Vatutin V. A. } (2012)  Limit theorems for weakly subcritical branching processes in random environment. Journal of Theoretical Probability, \textbf{25}, \textbf{  3}, 703--732.


\bibitem{AGKV2005} \textsc{ Afanasyev V.I., Geiger J.,  Kersting G. \&  Vatutin V.A.} (2005) Criticality for branching processes in random environment,  The Annals of Probability \textbf{ 33}, \textbf{2}, 645--673.
 
 

	\bibitem{BL}\textsc{Bougerol Ph.  \&   Lacroix J.} (1985) Products of Random Matrices with Applications to Schr\"odinger  Operators,  Birkh\"auser.
	


 
 \bibitem{DW}\textsc{ Denisov D.  \& Wachtel V.}  (2015)  Random walks in cones.  
Annals of Probability, {\bf 43},  992--1044.


\bibitem{DHKP}\textsc{Dolgopyat D.,  Hebbar P.,  Koralov L.  \&   Perlman M. } (2018)
Multi-type branching processes with time-dependent branching rates, Journal of Applied Probability, \textbf{  55},  \textbf{   3}, 701--727.


\bibitem{DGV}\textsc{Dyakonova E. E. and Geiger J.  \& Vatutin V. A.} (2004) On the survival probability and a functional limit theorem for
branching processes in random environment,  Markov Process and  Related Fields, {\bf 2}, 289--306 

\bibitem{DV} \textsc{Dyakonova E. E.  \& Vatutin V. A. } (2017) Multitype branching processes in random environment: survival probability for the critical case, Teor. Veroyatnost. i Primenen., 2017, {\bf  62},  {\bf  4},  634--653.

	
	\bibitem{FK}\textsc{Furstenberg H. \&  Kesten  H.}  (1960)  Product of random matrices.
Annals  Mathematical Statistics \textbf{31},  457--469. 

 

\bibitem{GK}\textsc{Geiger J.  \&   Kersting G.}  (2002) The survival  probability of a critical branching process in random environment.  Theory of Probability and Applications, \textbf{45},  \textbf{3}, 518--526.  
	 

	\bibitem{Hennion}\textsc{Hennion H.} (1997)  Limit theorems for products of positive random matrices, Annals of Probability     {\bf  25},    {\bf   4},  1545--1587.  
		
\bibitem{Jones}\textsc{Jones O. D. } (1997) On the convergence of multi-type branching processes with varying environments, Annals of Applied Probability, { \bf 7},    {\bf  3}, 772--801.
	
\bibitem{Kaplan}\textsc{Kaplan N.}  (1974) Some results about multidimensional branching processes with random environments, Annals of Probability,   \textbf{2},   441--455.

\bibitem{kersting2017}\textsc{Kersting G.} (2020) A unifying approach to branching processes in varying environments.   Journal of Applied Probability,   {\bf 57}, {\bf 1},  196--220.


	\bibitem{Ko}\textsc{Kozlov M. V. }  (1976) On the asymptotic behaviour of the probability of non-extinction for critical branching processes in a random environment. Theory of  Probability  and its Applications, {\bf    XXI},   \textbf{4},   791--804.
 
		
		\bibitem{LN} \textsc{Lamperti J.  \& Ney P.} (1968) Conditioned branching processes and their limiting diffusions,  Theory
Prob. Appl.   {\bf   13}, 128--139.

 


\bibitem{LPP}\textsc{Le Page E., Peign\'e M.  \& Pham C. }  (2018) The survival probability of a critical multi-type branching process in i.i.d. random environment,   Annals of Probability, 
  {\bf   46},   {\bf 5}, 2946--2972.

 	

\bibitem{PW} \textsc{Peign\'e M. \& Woess W. } (2019),   Recurrence of 2-dimensional queueing processes, and
random walk exit times from the quadrant,  preprint hal-02281986v1.


 	\bibitem{pham2018}\textsc{Pham C.} (2018)  Conditioned limit theorems for products of positive random matrices,   Lat. Am. J. Probab. Math. Stat.   {\bf   15}, 67--100.
	

 


\end{thebibliography}
\end{document}